\input amstex
\documentstyle{amsppt}
\nopagenumbers
\nologo
\magnification\magstephalf
\catcode`@=11
\redefine\output@{%
  \def\break{\penalty-\@M}\let\par\endgraf
  \global\voffset=-20pt
  \ifodd\pageno\global\hoffset=5pt\else\global\hoffset=-34pt\fi  
  \shipout\vbox{%
    \ifplain@
      \let\makeheadline\relax \let\makefootline\relax
    \else
      \iffirstpage@ \global\firstpage@false
        \let\rightheadline\frheadline
        \let\leftheadline\flheadline
      \else
        \ifrunheads@ 
        \else \let\makeheadline\relax
        \fi
      \fi
    \fi
    \makeheadline \pagebody \makefootline
  }%
  \advancepageno \ifnum\outputpenalty>-\@MM\else\dosupereject\fi
}
\font\cpr=cmr7
\newcount\xnumber
\footline={\xnumber=\pageno
\divide\xnumber by 7
\multiply\xnumber by -7
\advance\xnumber by\pageno
\ifnum\xnumber>0\hfil\else\vtop{\vskip 0.5cm
\noindent\cpr CopyRight \copyright\ Sharipov R.A.,
1995, 2006.}\hfil\fi}
\def\setfirstpage{\global\firstpage@true}
\catcode`\@=\active

\fontdimen3\tenrm=3pt
\fontdimen4\tenrm=0.7pt

\def\leaderfill{\leaders\hbox to 0.3em{\hss.\hss}\hfill}
\font\tvbf=cmbx12
\font\tvrm=cmr12
\font\etbf=cmbx8
\font\tencyr=wncyr10
\def\negskp{\hskip -2pt}

\def\MatGrGL{\operatorname{GL}}
\def\MatGrU{\operatorname{U}}
\def\MatGrSU{\operatorname{SU}}
\def\tr{\operatorname{tr}}
\def\End{\operatorname{End}}
\def\Hom{\operatorname{Hom}}
\def\Aut{\operatorname{Aut}}
\def\Cl{\operatorname{Cl}}
\def\Orbit{\operatorname{Orb}}
\def\Img{\operatorname{Im}}
\def\Ker{\operatorname{Ker}}
\def\compos{\,\raise 1pt\hbox{$\sssize\circ$} \,}
\Monograph
\loadbold
\TagsOnRight
\newcount\chapternum
\def\blue#1{#1}
\catcode`#=11\def\diez{#}\catcode`#=6
\catcode`_=11\catcode`_=8

\def\mytag#1{%
    \tag#1}
\def\mythetag#1{\thetag{\blue{#1}}\immediate\special{ps:
     ShrHPSdict begin /ShrBORDERthickness 0 def}}
\def\mythetagchapter#1#2{\thetag{\blue{#1}}\immediate\special{ps:
     ShrHPSdict begin /ShrBORDERthickness 0 def}}
\def\mybookref#1{#1}
\def\mybookcite#1{\cite{\blue{#1}}\immediate\special{ps:
     ShrHPSdict begin /ShrBORDERthickness 0 def}}
\def\myrefno#1{\no#1}
\def\myhref#1#2{\blue{#2}\immediate\special{ps:
     ShrHPSdict begin /ShrBORDERthickness 0 def}}
\def\myEarXivlink{\myhref{http://arXiv.org}{http:/\negskp/arXiv.org}}

\def\mytheorem#1{\csname proclaim\endcsname{Theorem #1}}
\def\mytheoremwithtitle#1#2{\csname proclaim\endcsname{Theorem #1#2}}
\def\mythetheorem#1{\blue{#1}\immediate\special{ps:
     ShrHPSdict begin /ShrBORDERthickness 0 def}}
\def\mythetheoremchapter#1#2{\blue{#1}\immediate\special{ps:
     ShrHPSdict begin /ShrBORDERthickness 0 def}}
\def\mylemma#1{\csname proclaim\endcsname{Lemma #1}}
\def\mylemmawithtitle#1#2{\csname proclaim\endcsname{Lemma #1#2}}
\def\mythelemma#1{\blue{#1}\immediate\special{ps:
     ShrHPSdict begin /ShrBORDERthickness 0 def}}
\def\mythelemmachapter#1#2{\blue{#1}\immediate\special{ps:
     ShrHPSdict begin /ShrBORDERthickness 0 def}}
\def\myproposition#1{\csname proclaim\endcsname{Proposition #1}}
\def\mypropositionwithtitle#1#2{\csname proclaim\endcsname{Proposition #1#2}}
\def\mytheproposition#1{\blue{#1}\immediate\special{ps:
     ShrHPSdict begin /ShrBORDERthickness 0 def}}

\def\mycorollary#1{\csname proclaim\endcsname{Corollary #1}}

\def\mydefinition#1{\definition{Definition #1}}
\def\mythedefinition#1{\blue{#1}\immediate\special{ps:
     ShrHPSdict begin /ShrBORDERthickness 0 def}}

\def\SectionNum#1#2{\S
     \,#1.}
\pagewidth{10cm}
\pageheight{15.5cm}
\fontdimen3\tenrm=3pt
\fontdimen4\tenrm=0.7pt

\document

\chapternum=1
\vbox to\vsize{\centerline{\etbf SPECIAL COURSE IN MATHEMATICS}
\vskip 3cm
\centerline{SHARIPOV\ R.\,A.}
\vskip 1.5cm
\centerline{\tvbf REPRESENTATIONS \ OF \ FINITE \ GROUPS}
\vskip 1.3cm
\centerline{\tvrm Part 1}
\vskip 6.5cm
\centerline{Ufa 1995}

\vss}
\newpage
\vbox to 13.5cm{
UDC 517.9\par
Sharipov R. A. {\bf Representations of finite groups. Part 1.} 
The textbook --- Ufa, 1995. --- 75 pages --- ISBN 5-87855-004-0.
\bigskip
\bigskip
     This book is an introduction to a fast developing branch of
mathematics --- the theory of representations of groups. It 
presents classical results of this theory concerning finite groups.
This book is written on the base of the special course which I gave 
in Mathematics Department of Bashkir State University.\par
     In preparing Russian edition of this book I used
computer typesetting on the base of \AmSTeX\ package and
I used cyrillic fonts of Lh-family distributed by CyrTUG
association of Cyrillic \TeX\ users. English edition is also
typeset by \AmSTeX.\par
\bigskip
\vfil
\line{ISBN 5-87855-004-0\hss\copyright\ Sharipov R.A., 1995}
\line{English Translation\hss\copyright\ Sharipov R.A., 2006}}
\newpage
\ \bigskip\medskip
\centerline{\bf CONTENTS.}
\bigskip
\line{CONTENTS.\ \leaderfill\ 3.}
\medskip
\line{PREFACE.\ \leaderfill\ \myhref{\diez pg4}{4}.}
\medskip
\line{CHAPTER~\uppercase\expandafter{\romannumeral 1}.
REPRESENTATIONS OF GROUPS.\ \leaderfill\ \myhref{\diez pg5}{5}.}
\medskip
\line{\S\,1. Representations of groups and their 
homomorphisms.\ \leaderfill\ \myhref{\diez pg5}{5}.}
\line{\S\,2. Finite dimensional representations.\ \leaderfill\ 
\myhref{\diez pg7}{7}.}
\line{\S\,3. Invariant subspaces. Restriction and factorization\hss}
\line{\qquad of representations.\ \leaderfill\ \myhref{\diez pg8}{8}.}
\line{\S\,4. Completely reducible representations.\ \leaderfill\ 
\myhref{\diez pg11}{11}.}
\line{\S\,5. Schur's lemma and some corollaries of it.\ \leaderfill\ 
\myhref{\diez pg21}{21}.}
\line{\S\,6. Irreducible representations of the direct product\hss}
\line{\qquad of groups.\ \leaderfill\ \myhref{\diez pg27}{27}.}
\line{\S\,7. Unitary representations.\ \leaderfill\ 
\myhref{\diez pg36}{36}.}
\bigskip
\line{CHAPTER~\uppercase\expandafter{\romannumeral 2}.
REPRESENTATIONS OF FINITE\hss}
 \line{\qquad GROUPS.\ \ \leaderfill\ \myhref{\diez pg44}{44}.}
\medskip
\line{\S\,1. Regular representations of finite groups.\ \leaderfill\ 
\myhref{\diez pg44}{44}.}
\line{\S\,2. Invariant averaging over a finite group.\ \leaderfill\ 
\myhref{\diez pg46}{46}.}
\line{\S\,3. Characters of group representations. \ \leaderfill\ 
\myhref{\diez pg50}{50}.}
\line{\S\,4. Orthogonality relationships.\ \leaderfill\ 
\myhref{\diez pg54}{54}.}
\line{\S\,5. Expansion into irreducible components.\ \leaderfill\ 
\myhref{\diez pg65}{65}.}
\bigskip
\line{CONTACTS.\ \ \leaderfill\ \myhref{\diez pg75}{75}.}
\bigskip
\line{APPENDIX.\ \ \leaderfill\ \myhref{\diez pg76}{76}.}
\newpage
\ \bigskip
\newtoks\truehead
\truehead=\headline
\headline{\hfill}
\centerline{\bf PREFACE.}
\bigskip
\medskip
     The theory of group representations is a wide branch of mathematics. 
In this book I explain very small part of this theory concerned with 
representations of finite groups. The material of the book is approximately
equal to a one semester course.\par
    In explaining the material of the book I tried to make it maximally
detailed, complete and self-consistent. The reader need not refer to other
literature. He is only required to know the linear algebra and the group
theory on the level of a standard university course. The linear algebra 
references are given to my book which now is available online:\par
\medskip
\noindent\vtop{\hsize=0.5cm\noindent [\mybookref{1}]}%
\vtop{\hsize=9.5cm\noindent Sharipov R. A. {\tencyr\char '074}Course of
linear algebra and multidimensional geometry{\tencyr\char '076}, Bashkir 
State University, Ufa 1996.}\par
\medskip
\noindent The primary material of the book is prepared on the base of 
the following two brilliant monographs:\par
\medskip
\noindent\vtop{\hsize=0.5cm\noindent [\mybookref{2}]}%
\vtop{\hsize=9.5cm\noindent
Neimark M. A. {\tencyr\char '074}Theory of representations of 
groups{\tencyr\char '076}, Nauka publishers, Moscow 1976;}\par
\medskip
\noindent\vtop{\hsize=0.5cm\noindent [\mybookref{3}]}%
\vtop{\hsize=9.5cm\noindent
Kirillov A. A. {\tencyr\char '074}Elements of the theory of 
representations{\tencyr\char '076}, Nauka publishers, Moscow 1978.}
\par
\medskip
Some problems concerning the character algebra for finite groups 
representations are not considered in this book. They will be given
in Part 2, which is planned as a separate book.
\bigskip\bigskip
\line{\vbox{\hsize 7.5cm\noindent September, 1995;\newline December,
2006.}\hss R.~A.~Sharipov.}
\newpage
\setfirstpage
\topmatter
\title\chapter{1}
REPRESENTATIONS OF GROUPS
\endtitle
\endtopmatter
\leftheadtext{CHAPTER \uppercase\expandafter{\romannumeral 1}.
REPRESENTATIONS OF GROUPS.}
\document
\headline=\truehead
\head
\SectionNum{1}{5} Representations of groups and their homomorphisms.
\endhead
\rightheadtext{\S\,1. Representations and their homomorphisms.}
     It is clear that matrix groups are in some sense simpler than abstract 
groups. Multiplication rule in them is explicitly specific. Dealing with
matrix groups one can use methods of the linear algebra and calculus. The
theory of group representations grows from the aim to reproduce an abstract
group in a matrix form.\par
     Let $V$ be a linear vector space over the field of complex numbers 
$\Bbb C$. By $\End(V)$ we denote the set of linear operators mapping $V$
into $V$. The subset of non-degenerate operators within $\End(V)$ is denoted
by $\Aut(V)$. It is easy to see that $\Aut(V)$ is a group. The operation of
composition, i\.\,e\. applying two operators successively, is the group
multiplication in $\Aut(V)$.
\mydefinition{1.1}A representation $f$ of a group $G$ in a linear vector
space $V$ is a group homomorphism $f\!:\,G\to\Aut(V)$.
\enddefinition
If $f$ is a representation of a group $G$ in $V$, this fact is briefly 
written as $(f,G,V)$. Let $g\in G$ be an element of the group $G$ , then
$f(g)$ is a non-degenerate operator acting within the space $V$. It is
called  the {\it representation operator\/} corresponding to the element 
$g\in G$. By $f(g)\bold x$ we denote the result of applying this operator 
to a vector $\bold x\in V$. The notation with two pairs of braces 
$f(g)(\bold x)$ will also be used provided it is more clear in a given
context. For instance, $f(g)(\bold x+\bold y)$. Representation operators 
satisfy the following evident relationships:
\roster
\item $f(g_1\,g_2)=f(g_1)\,f(g_2)$;
\item $f(1)=1$;
\item $f(g^{-1})=f(g)^{-1}$;
\endroster
\mydefinition{1.2} Let $(f,G,V)$ and $(h,G,W)$ be two representations
of the same group $G$. A linear mapping $A\!:\,V\to W$ is called a 
{\it homomorphism\/} sending the representation $(f,G,V)$ to the
representation $(h,G,W)$ if the following condition is fulfilled:
$$
\hskip -2em
A\compos f(g)=h(g)\compos A\text{\ \ for all \ }g\in G.
\mytag{1.1}
$$
The mapping $A$ in \mythetag{1.1}, which performs a homomorphism of two
representations, sometimes is called an {\it interlacing map}. 
\enddefinition
\mydefinition{1.3} A homomorphism $A$ interlacing two representations 
$(f,G,V)$ and $(h,G,W)$ is called an {\it isomorphism\/} if it is
bijective as a linear mapping $A\!:\,V\to W$.
\enddefinition
     It is easy to verify that the relation of being isomorphic is an
equivalence relation for representations. Two isomorphic representations
are also called {\it equivalent representations}. In the theory of
representations two isomorphic representations are treated as identical
representations because all essential properties of such two 
representations do coincide.
\mytheorem{1.1} If $(f,G,V)$ is a representation of a group $G$ in a
space $V$ and if $A\!:\,V\to W$ is a bijective linear mapping, then $A$
induces a unique representation of the group $G$ in $W$ which is equivalent
to $(f,G,V)$ and for which $A$ is an interlacing map.
\endproclaim
\demo{Proof} The proof of the theorem is trivial. Let's define the operators
of a representation $h$ in $W$ as follows:
$$
\hskip -2em
h(g)=A\compos f(g)\compos A^{-1}.
\mytag{1.2}
$$
It is easy to verify that the formula \mythetag{1.2} does actually define 
a representation of the group $G$ in $W$. Multiplying \mythetag{1.2} on 
the right by $A$, we get \mythetag{1.1}. Hence $A$ is an isomorphism 
of $f$ and $h$. Moreover, any representation of $G$ in $W$ for which $A$
is an interlacing isomorphism should coincide with $h$. This fact is proved
by multiplying \mythetag{1.1} on the right by $A^{-1}$.
\qed\enddemo
    From the theorem proved just above we conclude that if a representation
$f$ in $V$ is given, in order to construct an equivalent representation in
$W$ it is sufficient to have a linear bijection from $V$ to $W$. However,
in practice the problem is stated somewhat differently. Two representations
$f$ and $h$ in $V$ and $W$ are already given. The problem is to figure out
if they are equivalent and if so to find an interlacing operator. In this
statement this is one of the basic problems of the theory of representations.
\head
\SectionNum{2}{7} Finite dimensional representations.
\endhead
\rightheadtext{\S\,2. Finite-dimensional representations.}
\mydefinition{2.1} A representation $(f,G,V)$ is called 
finite-dimensional if its space $V$ is a finite-dimensional 
linear vector space, i\.\,e\. $\dim V<\infty$.
\enddefinition
     Below in this book we consider only finite-dimensional 
representations, though many facts proved for this case can then be 
transferred or generalized for the case of infinite-dimensional
representations.\par
     Note that any finite-dimensional linear vector space over the
field of complex numbers $\Bbb C$ can be bijectively mapped onto the
standard arithmetic coordinate vector space $\Bbb C^n$, where\linebreak 
$n=\dim V$. And $\Aut(\Bbb C^n)=\MatGrGL(n,\Bbb C)$. Therefore each
finite-dimensional representation is equivalent to some matrix
representation $f\!:\,G\to\MatGrGL(n,\Bbb C)$. This fact follows from 
the theorem~\mythetheorem{1.1}. In spite of this fact we shall consider
finite-dimensional representations in abstract vector spaces because
all statements in this case are more elegant and their proofs are sometimes
even more simple than the proofs of corresponding matrix statements.\par
\head
\SectionNum{3}{8} Invariant subspaces. Restriction and factorization of 
representations.
\endhead
\rightheadtext{\S\,3. Invariant subspaces.}
\mydefinition{3.1} Let $(f,G,V)$ be a representation of a group $G$ 
in a linear vector space $V$. A subspace $W\subseteq V$ is called an 
{\it invariant subspace\/} if for any $g\in G$ and for any $\bold x\in
W$ the result of applying $f(g)$ to $\bold x$ belongs to $W$, i\.\,e\.
$f(g)\bold x\in W$.
\enddefinition
    The concept of {\it irreducibility\/} is introduced in terms of invariant
subspaces. This is the central concept in the theory of representations.
\mydefinition{3.2} A representation $(f,G,V)$ of the group $G$ is called
{\it irreducible\/} if it has no invariant subspaces other than $W=\{0\}$ 
and $W=V$. Otherwise the representation $(f,G,V)$ is called {\it reducible}.
\enddefinition
    Assume that $(f,G,V)$ is an irreducible representation. Let's choose
some vector $\bold x\neq 0$ of $V$ and consider its orbit:
$$
\Orbit_f(\bold x)=\{\,\bold y\in V\!:\,\bold y=f(g)\bold x\text{\ \ for
some \ }g\in G\,\}.
$$
The orbit $\Orbit_f(\bold x)$ is a subset of the space $V$ invariant under
the action of the representation operators. However, in general case, it is
not a linear subspace. Let's consider its linear span
$$
W=\langle \Orbit_f(\bold x)\rangle.
$$
The subspace $W$ is invariant and $W\neq\{0\}$ since it possesses the 
non-zero vector $\bold x$. Then due to the irreducibility of $f$ we get
$W=V$. As a result one can formulate the following criterion of 
irreducibility.
\mytheoremwithtitle{3.1}{ ({\bf irreducibility criterion})} A representation
$(f,G,V)$ is irreducible if and only if the orbit of an arbitrary non-zero
vector $\bold x\in V$ spans the whole space $V$.
\endproclaim
    The necessity of this condition was proved above. Let's prove its 
sufficiency. Let $W\subseteq V$ be an invariant subspace such that 
$W\neq\{0\}$. Let's choose a nonzero vector $\bold x\in W$. Due to the
invariance of $W$ we have $\Orbit_f(\bold x)\subseteq W$. Hence, 
$\langle \Orbit_f(\bold x)\rangle\subseteq W$. But $\langle
\Orbit_f(\bold x)\rangle=V$, therefore $W=V$. The criterion is proved.
\par
    Irreducible representations are similar to chemical elements. One
cannot extract other more simple representations of a given group from
them. Any reducible representation in some sense splits into irreducible
ones. Therefore in the theory of representations the following two
problems are solved:
\rosteritemwd=0pt
\roster
\item to find and describe all irreducible representations of a given
group;
\item to suggest a method for splitting an arbitrary representation
into its irreducible components.
\endroster
The first problem is analogous to building the Mendeleev's table in
chemistry, the second problem is analogous to chemical analysis of
substances.\par
     Let's consider some reducible representation $(f,G,V)$ of a group
$G$. Assume that $W$ is an invariant subspace such that $\{0\}\subsetneq
W\subsetneq V$. Let's denote by $\varphi(g)$ the restriction of the
operator $f(g)$ to the subspace $W$:
$$
\hskip -2em
\varphi(g)=f(g)\,\hbox{\vrule height 8pt depth 10pt
width 0.5pt}_{\,W}.
\mytag{3.1}
$$
For the operators $\varphi(g)$ we have the following relationships:
$$
\align
&\hskip -2em\varphi(g)\,\varphi(g^{-1})=(f(g)\,f(g^{-1}))\,\hbox{\vrule 
height 8pt depth 10pt width 0.5pt}_{\,W}=1;
\mytag{3.2}\\
&\hskip -2em\varphi(g_1)\,\varphi(g_2)=(f(g_1)\,f(g_2))\,\hbox{\vrule 
height 8pt depth 10pt width 0.5pt}_{\,W}=\varphi(g_1\,g_2).
\mytag{3.3}
\endalign
$$
From the relationship \mythetag{3.2} we conclude that the operator 
$\varphi(g)$ is invertible and $\varphi(g)^{-1}=\varphi(g^{-1})$.
Hence, $\varphi(g)\in\Aut(W)$. The relationship \mythetag{3.3} in 
its turn shows that the mapping\linebreak $\varphi\!:\,G\to\Aut(W)$ 
is a group homomorphism defining a representation.\par
\mydefinition{3.3} The representation $(\varphi,G,W)$ of a group
$G$ obtained by restricting the operators of a representation 
$(f,G,V)$  to its invariant subspace $W\subseteq V$ according to
the formula \mythetag{3.1} is called the {\it restriction\/} of
$f$ to $W$.
\enddefinition
     The presence of an invariant subspace $W$ let's us define
factoroperators in the factorspace $V/W$:
$$
\hskip -2em
\psi(g)=f(g)\,\hbox{\vrule height 8pt depth 10pt
width 0.5pt}_{\,V/W}.
\mytag{3.4}
$$
Let's recall that the action of the operator $\psi(g)$ upon a coset
$\Cl_W(\bold x)$ in $V/W$ is defined as follows:
$$
\hskip -2em
\psi(g)\Cl_W(\bold x)=\Cl_W(f(g)\bold x).
\mytag{3.5}
$$
The correctness of the definition \mythetag{3.5} is verified by direct
calculations (see \mybookcite{1}). The factoroperators \mythetag{3.5} 
obey the following relationships:
$$
\align
&\hskip -2em
\gathered
\psi(g)\,\psi(g^{-1})\Cl_W(\bold x)=\Cl_W(f(g)f(g^{-1})\bold x)=\\
=\Cl_W(f(g\,g^{-1})\bold x)=\Cl_W(\bold x),
\endgathered
\mytag{3.6}\\
\vspace{2ex}
&\hskip -2em
\gathered
\psi(g_1)\,\psi(g_2)\Cl_W(\bold x)=\Cl_W(f(g_1)f(g_2)\bold x)=\\
=\Cl_W(f(g_1\,g_2)\bold x)=\psi(g_1\,g_2)\Cl_W(\bold x).
\endgathered
\mytag{3.7}
\endalign
$$
From \mythetag{3.6} and \mythetag{3.7} we conclude that the factoroperators 
\mythetag{3.4} satisfy the relationships similar to \mythetag{3.2} and
\mythetag{3.3}. They define a representation $(\psi,G,V/W)$ which is usually
called a {\it factorrepresentation}.\par
     The representations $(\varphi,G,W)$ and $(\psi,G,V/W)$ are generated 
by the representation $f$. Each of them inherits a part of the information 
contained in the representation $f$. In order to understand which part of 
the information is kept in $\varphi$ and $\psi$ let's study the
representation operators  $f(g)$ in some special basis. Let's choose a
basis $\bold e_1,\,\ldots,\,\bold e_s$ within the invariant subspace $W$.
Then we complete this basis up to a basis in the space $V$. This construction 
is based on the theorem on completing a basis of a subspace
(see \mybookcite{1}). The matrix of the operator $f(g)$ in such a composite 
basis is a block-triangular matrix:
$$
\hskip -2em
F(g)=\Vmatrix 
\hskip 0.2em \varphi^i_j\hskip 0.2em & \hskip -0.7em
\hbox{\vrule height 1.8ex depth 5.5ex}\hskip -0.4em & u^i_j\\
\vspace{-4.0ex}
\hskip -0.8em\vbox{\hsize 4em\hrule width 3.6em}\hskip -2.8em\\
0&\hskip -5em &\psi^i_j
\endVmatrix.
\mytag{3.8}
$$
The upper left diagonal block coincides with the matrix of the operator
$\varphi(g)$ in the basis $\bold e_1,\,\ldots,\,\bold e_s$. The lower right 
diagonal block coincides with the matrix of the factoroperator $\psi(g)$ in
the basis $\bold E_1,\,\ldots,\,\bold E_{n-s}$, where
$$
\bold E_1=\Cl_W(\bold e_{s+1}),\quad\bold E_2=\Cl_W(\bold e_{s+2}),
\quad.\ .\ .\ ,\quad\bold E_{n-s}=\Cl_W(\bold e_n).
$$
From \mythetag{3.8} we conclude that when passing from $f$ to its restriction 
$\varphi$ and to its factorrepresentation $\psi$ the amount
of the lost information is determined by the upper right non-diagonal
block $u^i_j$ in the matrix \mythetag{3.8}.\par
     Despite to the loss of information, the passage from $f$ to the
pair of representations $\varphi$ and $\psi$ can be treated as splitting
$f$ into more simple components. If $\varphi$ and $\psi$ are also
reducible, they can be split further. However, this process of splitting 
is finite since in each step we have the reduction of the dimension:
$\dim(W)<\dim(V)$ and $\dim (V/W)<\dim(V)$. The process will terminate
when we reach irreducible representations.\par
\head
\SectionNum{4}{11} Completely reducible representations.
\endhead
\rightheadtext{\S\,4. Completely reducible representations.}
    As we have seen above, the process of fragmentation of a reducible
representation leads to the loss of information. However, there is a
special class of representations for which the loss of information is
absent. This is the class of {\it completely reducible representations}.
\mydefinition{4.1} A representation $(f,G,V)$ of a group $G$ is called
{\it completely reducible\/} if each its invariant subspace $W$ has an
invariant direct complement $U$, i\.\,e\. $V$ is a direct sum of two
invariant subspaces $V=W\oplus U$.
\enddefinition
    Note that an irreducible representation is a trivial example of a
completely reducible one. Here wee have $V=V\oplus\{0\}$.\par
    Let $(f,G,V)$ be a reducible and completely reducible representation.
Let $W$ be an invariant subspace for $f$ and $U$ be its invariant direct
complement. Then we have the following isomorphisms of the restrictions
and factors:
$$
\xalignat 2
&\hskip -2em f\,\hbox{\vrule height 8pt depth 10pt
width 0.5pt}_{\,U}\cong f\,\hbox{\vrule height 8pt depth 10pt
width 0.5pt}_{\,V/W},
&f\,\hbox{\vrule height 8pt depth 10pt
width 0.5pt}_{\,W}\cong f\,\hbox{\vrule height 8pt depth 10pt
width 0.5pt}_{\,V/U}.
\mytag{4.1}
\endxalignat
$$
In order to prove \mythetag{4.1} lets consider again a basis $\bold e_1,
\,\ldots,\,\bold e_s$ in $U$ and complement it with a basis $\bold h_1,
\,\ldots,\,\bold h_{n-s}$ in $U$. Let's denote $\bold e_{s+1}=\bold h_1,
\,\ldots,\,\bold e_n=\bold h_{n-s}$. As a result of such a reunion of
bases we get a basis in $V$. Let's define a mapping $A\!:\,U\to V/W$ by
setting
$$
\hskip -2em
A\bold x=\Cl_W(\bold x)\text{\ \ for all \ }\bold x\in U.
\mytag{4.2}
$$
From \mythetag{4.2} one easily derives the values of the mapping $A$ on
basis vectors 
$$
\hskip -2em
A(\bold h_1)=\bold E_1,
\quad.\ .\ .\ ,\quad A(\bold h_{n-s})=\bold E_{n-s}.
\mytag{4.3}
$$
The mapping $A$ establishes bijective correspondence of bases of $U$ and
$V/W$. For this reason it is bijective. Let's verify that it implements
an isomorphism of representations, i\.\,e\. let's verify the relationship
\mythetag{1.1} for $A$:
$$
A(f(g)\bold h_i)=\Cl_W(f(g)\bold h_i)=\psi(g)\Cl_W(\bold h_i)=
\psi(g)A(\bold h_i).
$$
This relationship proves the first isomorphism in \mythetag{4.1}. It is
implemented by the mapping $A$ defined in \mythetag{4.2}.\par
     From the invariance of the subspace $U$ we derive $f(g)\bold h_i
\in U$. This fact and the relationship \mythetag{4.3} let us write the
matrix of the operator $f(g)$:
$$
\hskip -2em
F(g)=\Vmatrix 
\hskip 0.2em \varphi^i_j\hskip 0.2em & \hskip -0.7em
\hbox{\vrule height 1.8ex depth 5.5ex}\hskip -0.4em & 0\\
\vspace{-4.0ex}
\hskip -0.8em\vbox{\hsize 4em\hrule width 3.6em}\hskip -2.8em\\
0&\hskip -5em &\psi^i_j
\endVmatrix.
\mytag{4.4}
$$
The matrix \mythetag{4.4} is block-diagonal, therefore no loss of 
information occur when passing from $f$ to the representations 
$\varphi$ and $\psi$ in the case of completely reducible representation
$f$. Due to \mythetag{4.1} the representation $\psi$ can be treated
as the restriction of $f$ to the invariant complement $U$. The
representations $\varphi,G,W)$ and $(\psi,G,U)$ are not linked to each
other, they are defined in their own spaces which intersect trivially
and their sum is the complete space $V$. This situation is described by 
the following definition.
\mydefinition{4.2} A representation $(f,G,V)$ of a group $G$ is called 
an {\it inner direct sum\/} of the representations $(\varphi,G,W)$ and
$(\psi,G,W)$ if $V=W\oplus U$, while the subspaces $W$ and $U$ are
invariant subspaces for $f$ and the restrictions of $f$ to these 
subspaces coincide with $\varphi$ and $\psi$.
\enddefinition
     Note that the splitting of $f$ into a direct sum $f=\varphi\oplus
\psi$ can occur even in that case where $f$ is not a completely reducible 
representation. However in this case such a splitting is rather an
exception than a rule.\par
     Having a pair of representations $(\varphi,G,W)$ and $(\psi,G,U)$ 
of the same group $G$ in two distinct spaces, we can construct their 
{\it exterior direct sum}. Let's consider the external direct sum 
$W\oplus U$. Let's recall that it is the set of ordered pairs $(\bold
w,\bold u)$, where $\bold w\in W$ and $\bold u\in U$, with the algebraic
operations 
$$
\aligned
&(\bold w_1,\bold u_1)+(\bold w_2,\bold u_2)=(\bold w_1+\bold w_2,
\bold u_1+\bold u_2),\\
&\alpha\,(\bold w,\bold u)=(\alpha\,\bold w,\alpha\,\bold u)
\text{, \ \ where \ }\alpha\in\Bbb C.
\endaligned
$$
The subspaces $W$ and $U$ in the external direct sum $W\oplus U$ are
assumed to be disjoint even in that case where they have non-zero 
intersection or do coincide. Let's define operators $f(g)$ in 
$W\oplus U$ as follows:
$$
\hskip -2em
f(g)(\bold w,\bold u)=(\varphi(g)\bold w,\psi(g)\bold u).
\mytag{4.5}
$$
The representation $(f,G,W\oplus U)$ constructed according to the 
representation \mythetag{4.5} is called the {\it external direct sum\/}
of the representations $(\varphi,G,W)$ and $(\psi,G,W)$. It is denoted
$f=\varphi\oplus\psi$. The difference of internal and external direct
sums is rather formal, their properties do coincide in most.\par
     Assume that the space $V$ of a representation $(f,G,V)$ is expanded
into a direct sum of subspaces $V=W\oplus U$ (not necessarily invariant
ones). Each such an expansion are uniquely associated with two projection
operators $P$ and $Q$. The projector $P$ is the projection operator that
projects onto the subspace $W$ parallel to the subspace $U$, while $Q$
projects onto $U$ parallel to $W$. They satisfy the following 
relationships:
$$
\xalignat 3
&\hskip -2em
P^2=P,
&&Q^2=Q,
&&P+Q=1.
\quad
\mytag{4.6}
\endxalignat
$$
Moreover, $W=\Img P$ and $U=\Img Q$. These properties of the projection
operators are well known (see \mybookcite{1}).\par
\mylemma{4.1} The subspace $W$ in an expansion $V=W\oplus U$ is invariant
with respect to the operators of a representation $f,G,V)$ if and only if
the corresponding projector $P$ obeys the relationship
$$
\hskip -2em
(P\compos f(g)-f(g)\compos P)\compos P=0\text{\ \ for all \ }g\in G.
\mytag{4.7}
$$
The invariance of both subspaces $W$ and $U$ in an expansion $V=W\oplus U$
is equivalent to the relationship 
$$
\hskip -2em
(P\compos f(g)=f(g)\compos P)\text{\ \ for all \ }g\in G,
\mytag{4.8}
$$
which means that the projector $P$ commutes with all operators of the
representation $f$.
\endproclaim
    Let $\bold x$ be an arbitrary vector. Then $P\bold x\in W$. In the 
case of an invariant subspace $W$ the vector $\bold y=f(g)P\bold x$ also
belongs to $W$. For the vector $\bold y\in W$ we have $P\bold y=\bold y$.
Therefore
$$
\hskip -2em
Pf(g)P\bold x=f(g)P\bold x=f(g)P^2\bold x.
\mytag{4.9}
$$
Comparing the left and right hand sides of the equality \mythetag{4.9}
and taking into account the arbitrariness of the vector $\bold x$, we
easily derive the relationship \mythetag{4.7} in the statement of the 
lemma.\par
     And conversely, from the relationship \mythetag{4.7}, using the
property $P^2=P$ from \mythetag{4.6}, we easily derive \mythetag{4.9}.
From \mythetag{4.9} we derive that the vector $f(g)P\bold x$ belongs
to $W$ for an arbitrary vector $\bold x$. Let's choose $\bold x\in W$
and from $P\bold x=\bold x$ for such a vector $\bold x$ we find that
$f(g)\bold x$ belongs to $W$. The invariance of $W$ is established. In
order to prove the second statement of the lemma let's write the
relationship \mythetag{4.7} for the projector $Q$:
$$
\hskip -2em
(Q\compos f(g)-f(g)\compos Q)\compos Q=0.
\mytag{4.10}
$$
from $Q=1-P$ we get $Q\compos f(g)-f(g)\compos Q=f(g)\compos P-P
\compos f(g)$. Therefore the relationship \mythetag{4.10} is rewritten
as
$$
f(g)\compos P-P\compos f(g)+(P\compos f(g)-f(g)\compos P)\compos P=0.
$$
Then, taking into account \mythetag{4.7}, we reduce it to the relationship
\mythetag{4.8}, which means that the operators $P$ and $f(g)$ do commute.
The lemma is proved.\par
     The second proposition of the lemma~\mythelemma{4.1} can be 
generalized for the case where $V$ is expanded into a direct sum of
several subspaces. Assume that $V=W_1\oplus\,\ldots\,\oplus W_s$.
Let's recall (see \mybookcite{1}) that this expansion uniquely fixes 
a concordant family of projection operators $P_1,\,\ldots,\,P_s$. 
They obey the following concordance relationships:
$$
\xalignat 2
&(P_i)^2=P_i,&& P_1+\ldots+P_s=1,\\
&P_i\compos P_j=0\text{\ \ for \ }i\neq j,
&&P_i\compos P_j=P_j\compos P_i.
\endxalignat
$$
Moreover, $W_i=\Img P_i$. The condition of invariance of all subspaces
$W_i$ in the expansion $V=W_1\oplus\,\ldots\,\oplus W_s$ with respect to
representation operators of a representation $(f,G,V)$ in terms of the
corresponding projection operators is formulated in the following lemma.
We leave its proof to the reader.
\mylemma{4.2} The expansion $V=W_1\oplus\ldots\oplus W_s$ is an expansion
of $V$ into a direct sum of invariant subspaces of the representation
$(f,G,V)$ if and only if all projection operators $P_i$ associated with the
expansion $V=W_1\oplus\ldots\oplus W_s$ commute with the representation
operators $f(g)$.
\endproclaim
    Let's consider a completely reducible finite-dimensional representation
$(f,G,V)$ split into a direct sum of its restrictions to invariant subspaces
$f=\varphi\oplus\psi$. Assume that the restriction $\varphi$ of $f$ to the
subspace $W$ is reducible and assume that $W_1$ is its non-trivial invariant
subspace: $\{0\}\subsetneq W_1\subsetneq W$. Then the subspace $W_1$ is
invariant with respect to $f$. It has an invariant complement $U_1$. Let's
consider the expansions
$$
\xalignat 2
&\hskip -2em
V=W\oplus U,
&&V=W_1\oplus U_1.
\mytag{4.11}
\endxalignat
$$
Note that $W_1\subset W$, hence, $W+U_1=V$. Therefore the dimensions of the 
subspaces in \mythetag{4.11} are related as follows:
$$
\dim(W)+\dim(U_1)-\dim(V)=\dim(W)-\dim(W_1)>0.
$$
From this relationship, we see that the intersection $W_2=W\cap U_1$ is
non-zero and its dimension is given by the formula
$$
\hskip -2em
\dim(W\cap U_1)=\dim(W)-\dim(W_1).
\mytag{4.12}
$$
Now note that $W_1\cap W_2=\{0\}$ since $W_2\subset U_1$. Therefore, due to
\mythetag{4.12} the subspace $W$ is expanded into the direct sum
$$
\pagebreak
W=W_1\oplus W_2,
$$
each summand in which is invariant with respect to $f$ and, hence, with 
respect to $\varphi$. Thus, we have proved the following important 
theorem.
\mytheorem{4.1} The restriction of a completely reducible finite-dimensional
representation to an invariant subspace is completely reducible.
\endproclaim
The next theorem on the expansion into a direct sum is an immediate 
consequence of the theorem~\mythetheorem{4.1}.
\mytheorem{4.2} Each finite-dimensional completely reducible representation
$f$ is expanded into a direct sum of irreducible representations 
$$
\xalignat 2
&\hskip -2em
f=f_1\oplus\ldots\oplus f_k,
&&V=W_1\oplus\ldots\oplus W_k,
\quad
\mytag{4.13}
\endxalignat
$$
where each $f_i$ is a restriction of $f$ to the corresponding invariant 
subspace $W_i$.
\endproclaim
Note that in general the expansion \mythetag{4.13} is not unique. Let's
consider two expansions of $f$ into irreducible components:
$$
\xalignat 2
&\hskip -2em
f=f_1\oplus\ldots\oplus f_k,
&&V=W_1\oplus\ldots\oplus W_k,
\quad\\
\vspace{-1.2ex}
&&&\mytag{4.14}\\
\vspace{-1.2ex}
&\hskip -2em
f=\tilde f_1\oplus\ldots\oplus\tilde f_k,
&&V=\tilde W_1\oplus\ldots\oplus\tilde W_k.
\quad
\endxalignat
$$
The extent of differences in two expansions \mythetag{4.14} is determined
by the following Jordan-H\"older theorem.
\mytheorem{4.3} The numbers of irreducible components in the expansions
\mythetag{4.14} are the same: $q=k$, and there is a transposition 
$\sigma \in S_k$ such that $(f_i,G,W_i)\cong (\tilde f_{\sigma i},G,
\tilde W_{\sigma i})$.
\endproclaim
     The expansions \mythetag{4.14} are isomorphic up to a transposition
of components. However, we should emphasize that the isomorphism does
not mean the coincidence of these expansions.\par
     We shall prove the Jordan-H\"older theorem by induction on the number
of components $k$ in the first expansion \mythetag{4.14}.\par
     {\bf The base of the induction}: $k=1$, $V=W_1$. In this case the
representation $f=f_1$ is irreducible. Therefore $q=1=k$ and $\tilde W_1
=V$, $\tilde f_1=f=f_1$. The base of the induction is proved.\par
     {\bf The inductive step}. Assume that the theorem is valid for
representations possessing at least one expansion of the form
\mythetag{4.13} with the length $k-1$. For the representation $f$ we
introduce the following notations:
$$
\hskip -2em
\aligned
&\tilde V_i=\tilde W_1\oplus\ldots\oplus\tilde W_i\text{\ \ where \ }
i=1,\,\ldots,\,q,\\
&U=W_1\oplus\ldots\oplus W_{k-1}\text{\ \ and \ }\tilde U_i=\tilde V_i
\cap U.
\endaligned
\mytag{4.15}
$$
All of the subspaces in \mythetag{4.15} are invariant with respect to $f$.
Moreover, $V=U\oplus W_k$ and there are two chains of inclusions:
$$
\hskip -2em
\aligned
&\{0\}\subsetneq\tilde V_1\subsetneq\ldots\subsetneq\tilde V_q=V,\\
&\{0\}\subsetneq\tilde U_1\subsetneq\ldots\subsetneq\tilde U_q=U.
\endaligned
\mytag{4.16}
$$
Let $h\!:\,V\to V/U$ be a canonical projection onto the factorspace. Let's
denote by $\varphi$ the factorrepresentation in the factorspace $V/U$ and
consider the chain of invariant subspaces for it
$$
\{0\}\subseteq h(\tilde V_1)\subseteq\ldots\subseteq h(\tilde V_q)
=h(V)=V/U\cong W_k.
$$
Due to the isomorphism $\varphi\cong f_k$ we conclude that $\varphi$ is
irreducible. Hence the above chain does actually look like 
$$
\{0\}=h(\tilde V_1)=\ldots=h(\tilde V_s)\subsetneq h(\tilde V_{s+1})
=\ldots=h(\tilde V_q)=h(V).
$$
Therefore $\tilde V_i\subseteq U$ and $\tilde U_i=\tilde V_i$ for
$i\leqslant s$. For $i\geqslant s+1$ we use the isomorphisms 
$\tilde V_i/\tilde U_i\cong h(\tilde V_i)=h(\tilde V_{i+1})\cong
\tilde V_{i+1}/\tilde U_{i+1}$. But $\tilde V_i\subsetneq\tilde V_i$, 
hence, $\tilde U_i\subsetneq\tilde U_{i+1}$. Then from \mythetag{4.16} 
we get
$$
\hskip -2em
\{0\}\subsetneq\tilde U_1\subsetneq\ldots\subsetneq\tilde U_s
=\tilde U_{s+1}\subsetneq\ldots\subsetneq\tilde U_q=U.
\mytag{4.17}
$$
The equality $\tilde U_s=\tilde U_{s+1}$ follows from the irreducibility
of the factorrepresentation $\varphi_{s+1}\cong\varphi\cong f_k$ in the
factorspace $\tilde V_{s+1}/\tilde V_s$ and from the inclusions $\tilde 
V_s=\tilde U_s=\tilde U_{s+1}\subsetneq\tilde V_{s+1}$. Relying upon 
complete reducibility of $f$, we choose invariant complements $\tilde
W_{i+1}$ for $\tilde U_i$ in $\tilde U_{i+1}$, i\.\,e\. $\tilde U_{i+1}
=\tilde U_i\oplus\tilde W_{i+1}$. Then $\tilde V_i+\tilde U_{i+1}=\tilde
V_i\oplus\tilde W_{i+1}$ is an invariant subspace in $\tilde V_{i+1}$
containing $\tilde V_i$ and not coinciding with it. But the 
factorrepresentation in $\tilde V_{i+1}/\tilde V_i$ is isomorphic to
$tilde f_{i+1}$ and, hence, is  irreducible. Therefore $\tilde V_i
\oplus\tilde W_{i+1}=\tilde V_{i+1}$ and the restriction of $f$ to
$\tilde W_{i+1}$ is isomorphic to $\tilde f_{i+1}$. From \mythetag{4.15} 
and \mythetag{4.17} we get
$$
\hskip -2em
\aligned
&U=W_1\oplus\ldots\oplus W_s\oplus W_{s+1}\oplus\ldots\oplus W_{k-1},\\
&U=\tilde W_1\oplus\ldots\oplus\tilde W_s\oplus\tilde W_{s+2}\oplus
\ldots\oplus\tilde W_q.
\endaligned
\mytag{4.18}
$$
Now it is sufficient to apply the inductive hypothesis to \mythetag{4.18}.
As a result we get $k=q$ and find $\sigma i$ for $i=1,\,\ldots,\,k-1$. The
isomorphism of $f_k$ and the factorrepresentation $\varphi_{s+1}$ in
$\tilde V_{s+1}/\tilde V_s$ yields $\sigma k=s+1$ since $\varphi_{s+1}
\cong\tilde f_{s+1}$ by construction of the subspaces \mythetag{4.15}.
The Jordan-H\"older theorem is proved.\par
\smallskip
     Due to the theorems~\mythetheorem{4.1} and \mythetheorem{4.2} it is
natural to introduce the following terminology.\par
\mydefinition{4.3} An invariant subspace $W\subseteq V$ is called an
{\it irreducible subspace\/} for a representation $(f,G,V)$ if the
restriction of $f$ to $W$ is irreducible.
\enddefinition
      The following theorem yields a tool for verifying the complete
reducibility of representations.\par
\mytheorem{4.4} A finite-dimensional representation $(f,G,V)$ is
completely reducible if and only if the set of all its irreducible
subspaces span the whole space $V$.
\endproclaim
\demo{Proof} Let $\{W_\alpha\}_{\alpha\in A}$ be the set of all irreducible
subspaces in $V$. The number of such subspaces could be infinite even in the
case of a finite-dimensional representation. \pagebreak Let's denote by $W$
the sum of all irreducible subspaces $W_\alpha$:
$$
W=\sum_{\alpha\in A}W_\alpha=\lower 3pt\hbox{$\left<\raise 3pt
\hbox{$\dsize\bigcup_{\alpha\in A}W_\alpha$}\right>$}.
$$
The theorem~\mythetheorem{4.4} says that the condition $W=V$ is necessary
and sufficient for the representation $f$ to be completely reducible. The
necessity of this condition follows from the theorem~\mythetheorem{4.2}.
Let's prove its sufficiency. Let $U=U_0$ be an invariant subspace for $f$
and $\{0\}\neq U_0\neq V$. The for each irreducible subspace $W_\alpha$
we have two mutually exclusive options:
$$
W_\alpha\subseteq U_0\text{\ \ or \ }W_\alpha\cap U_0=\{0\}.
$$
And there is at least one subspace $W_{\alpha_1}$ for which the second 
option is valid: $W_{\alpha_1}\cap U_0=\{0\}$. Indeed, otherwise we would
have $W\subseteq U_0$, which contradict $W=V$. Let's consider the sum
$$
U_1=U_0+W_{\alpha_1}=U_0\oplus W_{\alpha_1}.
$$
The subspace $U_1$ is invariant with respect to $f$. If $U_1\neq V$, we
repeat our considerations with $U_1$ instead of $U_0$. As a result we get
a new invariant subspace
$$
U_2=U_1\oplus W_{\alpha_2}=U_0\oplus W_{\alpha_1}\oplus W_{\alpha_2}.
$$
The process of adding new direct summand to $U_0$ will terminate at some
step since $\dim V<\infty$. As a result we get
$$
U_k=U_0\oplus(W_{\alpha_1}\oplus\ldots\oplus W_{\alpha_k})=V.
$$
The the subspace $W=W_{\alpha_1}\oplus\ldots\oplus W_{\alpha_k}$ is a
required invariant direct complement for $U=U_0$. The theorem is proved.
\qed\enddemo
\head
\SectionNum{5}{21} Schur's lemma and some corollaries of it.
\endhead
\rightheadtext{\S\,5. Schur's lemma and some corollaries \dots}
    The theorem~\mythetheorem{4.2} proved in previous section says that
each completely reducible representation is expanded into a direct sum of
its  irreducible components. In this section we study these irreducible
components by themselves. Schur's lemma plays the central role in
this. We formulate two versions of this lemma, the second version being 
a strengthening of the first one, but for a more special case.
\mylemmawithtitle{5.1}{ ({\bf Schur's lemma})} Let $(f,G,V)$ and $(h,G,W)$
be two irreducible representations of a group $G$. Each homomorphism $A$
relating these two representations either is identically zero or is a
homomorphism.
\endproclaim
    Before proving the Schur's lemma we consider the following theorem
having a separate value.
\mytheorem{5.1} If a linear mapping $A\!:\,V\to W$ is a homomorphism of
representations from $(f,G,V)$ to $(h,G,W)$, then its kernel $Ker A$ is
an invariant subspace for $f$, while its image $\Img A$ is an invariant
subspace for $h$.
\endproclaim
\demo{Proof} Let $\bold y=A\bold X$ be a vector from the image of the
mapping $A$. We apple the operator $f(g)$ to it and use the relationship
\mythetag{1.1}. Then we get
$$
h(g)\bold y=h(g)A\bold x=Af(g)\bold x.
$$
Now it is clear that $h(g)\bold y\in\Img A$, i\.\,e\. $\Img A$ is an
invariant subspace for $h$. The invariance of $\Ker A$ is proved similarly.
Assume that $\bold x\in\Ker A$. Let's apply the operator $f(g)$ to $\bold x$
and then apply the mapping $A$. Taking into account \mythetag{1.1}, we get
$$
Af(g)\bold x=h(g)A\bold x=0
$$
since $A\bold x=0$. Hence, $f(g)\bold x\in\Ker A$. The invariance of the
kernel $\Ker A$ with respect to $f$ is proved.
\qed\enddemo
    Now let's proceed to proving Schur's lemma~\mythelemma{5.1}. The case, 
where $A=0$, is trivial. In order to exclude this case assume that
$A\neq 0$. Then $\Img A\neq\{0\}$. According to the
theorem~\mythetheorem{5.1} proved just above, $\Img A$ is invariant with
respect to the representation $h$. Since $h$ is irreducible, we get
$\Img A=W$. This means that the homomorphism $A$ is a surjective linear 
mapping $A\!:\,V\to W$.\par
    The kernel $\Ker A$ of the mapping $A\!:\,V\to W$ is invariant with
respect to $f$. Since $f$ is an irreducible representation, we have two
options $\Ker A=\{0\}$ or $\Ker A=V$. The second option leads to the 
trivial case $A=0$, which is excluded. Therefore, $\Ker A=\{0\}$. This
means that $A\!:\,V\to W$ is an injective linear mapping. Being surjective
and injective simultaneously, this mapping is bijective. Hence, it is 
an isomorphism of $(f,G,V)$ and $(h,G,W)$. Schur's lemma is proved.\par
    Now let's consider an operator $A\!:\,V\to V$ that interlaces an
irreducible representation $(f,G,V)$ with itself. The relationship
\mythetag{1.1} written for this case means that $A$ commute with all
representation operators $f(g)$. The second version of Schur's lemma
describes this special case $f=g$.
\mylemmawithtitle{5.2}{ ({\bf Schur's lemma})} Each operator $A\!:\,V
\to V$ that commutes with all operators of an irreducible 
finite-dimen\-sional representation $(f,G,V)$ in a linear vector 
space $V$ over the field of complex numbers $\Bbb C$ is a scalar 
operator, i\.\,e\. $A=\lambda\cdot 1$, where $\lambda\in\Bbb C$.
\endproclaim
\demo{Proof} Let $\lambda$ be an eigenvalue of the operator $A$ and
let $V_\lambda\neq\{0\}$ be the corresponding eigenspace. If $\bold x
\in V_\lambda$ then $A\bold x=\bold x$. Moreover, from the commutation
condition of the operators $A$ and $f(g)$ we obtain
$$
Af(g)\bold x=f(g)A\bold x=\lambda f(g)\bold x.
$$
Hence, $f(g)\bold x$ is also a vector of the eigenspace $V_\lambda$. In
other words, $V_\lambda$ is an invariant subspace. Since $V_\lambda\neq
\{0\}$ and since $f$ is irreducible, we get $V_\lambda=V$. Hence, 
$A\bold x=\bold x$ for an arbitrary vector $\bold x\in V$. This means 
that $A=\lambda\cdot 1$.
\qed\enddemo
    Note that the condition of finite-dimensionality of the representation
$f$ and the condition of complexity of the vector space $V$ in Schur's 
lemma~\mythelemma{5.2} are essential. Without these conditions one cannot
grant the existence of a non-trivial eigenspace.\par
    Now we apply Schur's lemma for to investigate tensor products of
representations of some special sort. For this reason we need to define
the concept of tensor product for representations.\par
    Let $(f,g,V)$ and $(h,G,W)$ be two representations of a group $G$. We
define the operators $\varphi(g)$ acting within the tensor product
$V\otimes W$ by setting
$$
\hskip -2em
\varphi(g)(\bold v\otimes\bold w)=(f(g)\bold v\otimes h(g)\bold w)
\text{\ \ for all \ }g\in G.
\mytag{5.1}
$$
The operators $\varphi(g)$ constitute a new representation of the group $G$.
The proof of this proposition and verifying the correctness of the formula
\mythetag{5.1} are left to the reader.
\mydefinition{5.1} The representation $(\varphi,G,V\otimes W)$ given by 
the operators \mythetag{5.1} is called the {\it tensor product\/} of the
representations $f$ and $h$. It is denoted $\varphi=f\otimes h$.
\enddefinition
     Note that the construction \mythetag{5.1} is easily generalized for
the case of several representations.\par
     Assume that the representation $f$ in the construction \mythetag{5.1}
is irreducible. As a second tensorial multiplicand in \mythetag{5.1} we
choose the trivial representation $i$ given by the identity operators
$i(g)=1$ for all $g\in G$. In the case where $\dim W>1$ the tensor product
$f\otimes i$ is a reducible representation. Let's prove this fact. Assume 
that $\bold e_1,\,\ldots,\,\bold e_m$ is a basis in the space $W$. Let's
denote by $W_k$ one-dimensional subspaces spanned by basis vectors 
$\bold e_k$ and then consider the subspaces
$$
\hskip -2em
V\otimes W_k,\quad k=1,\,\ldots,\,m.
\mytag{5.2}
$$
The subspaces \mythetag{5.2} are invariant within $V\otimes W$ with 
respect to the operators of the representation $f\otimes i$.The 
restrictions of $f\otimes i$ to $V\otimes W_k$ all are isomorphic to
$f$. Hence, the subspaces \mythetag{5.2} are irreducible. These results 
are simple. They are proved immediately.\par
     No note that $V\otimes W=V\otimes W_1\oplus\ldots\oplus V\otimes W_m$.
The space $V\otimes W$ is a sum of its irreducible subspaces $V\otimes W_k$.
Hence it is spanned by these subspaces. Therefore, it is sufficient to
apply the theorem~\mythetheorem{4.4}. As a result we get the following 
proposition.
\mytheorem{5.2} The tensor product $f\otimes i$ of a finite-dimensional
irreducible representation $f$ and a trivial finite-dimensional
representation $i$ is completely reducible.
\endproclaim
By means of this theorem one can describe the structure of all invariant 
subspaces of the tensor product $f\otimes i$.
\mytheorem{5.3} If under the assumptions of the theorem~\mythetheorem{5.2}
$f$ and $i$ are representations in complex linear vector spaces $V$ and $W$, 
then each invariant subspace of the representation $f\otimes i$ is
a tensor product $U=V\otimes \tilde W$, where $\tilde W$ is some subspace
of $W$.
\endproclaim
\demo{Proof} Let $U$ be some invariant subspace of $f\otimes i$ in 
$V\otimes W$. Due to the theorem~\mythetheorem{5.2} the representation
$f\otimes i$ is completely reducible. Therefore $U$ has an invariant
complement $\tilde U$. The expansion $V\times W=U\oplus\tilde U$ means 
that we can define a projection operator $P$ that project onto $U$ 
parallel to $\tilde U$. The invariance of both spaces $U$ and $\tilde U$
means that $P$ commutes with all operators of the representation
$\varphi=f\times i$.\par 
     Let $A\!:\,V\times W\to V\times W$ be some operator in the tensor 
product $V\times W$. In our proof $A=P$, however, we return to this
special case a little bit later. Let's apply the operator $A$ to 
$\bold x\otimes\bold e_j$ and write the result as
$$
\hskip -2em
A(\bold x\otimes\bold e_j)=\sum^m_{j=1}A^k_j(\bold x)\otimes\bold e_k.
\mytag{5.3}
$$
Here $\bold e_1,\,\ldots,\,\bold e_m$ is some basis in $W$. Since basis
vectors are linearly independent, the coefficients $A^k_j(\bold x)\in V$
in the expansion \mythetag{5.3} are unique. They are changed only if we 
change a basis, the transformation rule for them being analogous to that
for components of a tensor. This fact is not important for us, since in
our considerations we will not change the basis $\bold e_1,\,\ldots,\,
\bold e_m$.\par
     Let's study the dependence of $A^k_j(\bold x)$ on $\bold x$. It is
clear that it is linear. Therefore each coefficient $A^k_j(\bold x)$ in
\mythetag{5.3} defines a linear operator $A^k_j\!:\,V\to V$. Let's write 
the commutation relationship for $A$ and $\varphi(g)$. It is sufficient
write this relationship as applied to the vectors of the form
$\bold x\otimes\bold e_i$. From \mythetag{5.3} we derive
$$
\hskip -2em
\aligned
&A\compos\varphi(g)(\bold x\otimes\bold e_i)=\sum^m_{k=1}
A^k_j\compos f(g)(\bold x)\otimes \bold e_k,\\
&\varphi(g)\compos A(\bold x\otimes\bold e_i)=\sum^m_{k=1}
f(g)\compos A^k_j(\bold x)\otimes \bold e_k.
\endaligned
\mytag{5.4}
$$
For to write \mythetag{5.4} it is sufficient to recall that $\varphi=
f\otimes i$, where $i$ is a trivial representation. From \mythetag{5.4},
it is easy to see that the commutation relationship $A\compos\varphi(g)
=\varphi(g)\compos A$ is equivalent to 
$$
\hskip -2em
f(g)\compos A^k_j=A^k_j\compos f(g).
\mytag{5.5}
$$
Thus, all of the linear operators $A^k_j$ commute with the operators
$f(g)$ of the representation $f$.\par
     The next step in proof is to return to the projection operator $A=P$
(see above) and to apply Schur's lemma to the operators $A^k_j=P^k_j$. The
projector $P$ commutes with $\varphi(g)$, therefore the relationships
\mythetag{5.5} hold for $A^k_j=P^k_j$. Applying Schur's
lemma~\mythelemma{5.2}, we get
$$
\hskip -2em
P^k_j=\lambda^k_j\times 1,
\mytag{5.6}
$$
where $\lambda^k_j$ are some complex numbers. Substituting \mythetag{5.6}
into the expansion \mythetag{5.3}, for the operator $P$ we derive:
$$
\hskip -2em
P(\bold x\otimes\bold e_j)=\sum^m_{k=1}A^k_j(\bold x)\otimes\bold e_k=
\bold x\otimes\left(\,\shave{\sum^m_{k=1}}\lambda^k_j\,\bold e_k\right)\!.
\mytag{5.7}
$$
The relationship \mythetag{5.7} shows that it is natural to define a linear 
operator $Q\!:\,W\to W$ given by its values on basis vectors:
$$
\hskip -2em
Q(\bold e_j)=\sum^m_{k=1}\lambda^k_j\,\bold e_k.
\mytag{5.8}
$$
Due to \mythetag{5.8} we can rewrite \mythetag{5.7} as
$$
\hskip -2em
P(\bold x\otimes\bold e_j)=\bold x\otimes Q(\bold e_j).
\mytag{5.9}
$$
\par
Now let's remember that $P$ is a projection operator. Hence, $P^2=P$
(see \mybookcite{1}). Combining this relationship with \mythetag{5.9}, 
we derive $Q^2=Q$. Therefore, $Q$ is also a projection operator. Let's 
denote $\tilde W=\Img Q\subseteq W$. Then for the invariant subspace
$U\subseteq V\otimes W$ we get
$$
\hskip -2em
U=\Img P=V\otimes\Img Q=V\otimes\tilde W.
\mytag{5.10}
$$
The formula \mythetag{5.10} describes the structure of all invariant
subspaces if the representation $f\otimes i$. Thus, the 
theorem~\mythetheorem{5.3} is proved.
\qed\enddemo
    The theorem~\mythetheorem{5.3} can be applied for proving the
following extremely useful fact.
\mytheorem{5.4} Let $f$ be a finite dimensional irreducible representation
of some group $G$ in a complex linear vector space $V$. Then the set of
representation operators $f(g)$ spans the whole space of linear operators
$\End(V)$.
\endproclaim
\demo{Proof} Each representation $f$ in $V$ generates an associated
representation in the space of linear operators $\End(V)$. Indeed, if
$A\in\End(V)$, then we can define the action of $\psi(g)$ to $A$ as
the composition of operators $f(g)$ and $A$:
$$
\hskip -2em
\psi(g)(A)=f(g)\compos A\text{\ \ for all \ }A\in\End(V).
\mytag{5.11}
$$
Let $F$ be the span of the set of all operators $f(g)$, i\.\,e\.
$$
F=\left<\{f(g)\!:\,g\in G\}\right>. 
$$
It is easy to verify that the subspace $F$ is invariant with respect 
to the representation $\psi$ defined by the formula \mythetag{5.11}.
\par
     In order to apply the theorem~\mythetheorem{5.3} let's remember
that there is a canonical isomorphism $V\otimes V^*\cong\End(V)$,
where $V^*$ is the dual space for $V$ (the space of linear functionals
in $V$). This isomorphism is established by the mapping 
$\sigma\!:\,V\times V^*\to\End(V)$ which is defined as follows:
$$
\sigma(\bold x\otimes\lambda)\bold y=\lambda(y)\bold x
\text{\ \ for all \ }\bold x,\bold y\in V\text{ \ \ and \ }
\lambda\in V^*.\quad
\mytag{5.12}
$$
The proof of the correctness of the definition \mythetag{5.12} and 
verifying that $\sigma$ is an isomorphism are left to the reader. 
\par
    It is easy to verify that the canonical isomorphism $\sigma$ 
is an isomorphism interlacing the representation $\psi$ from 
\mythetag{5.11} and the representation $f\times i$, where $i$ is
the trivial representation of the group $G$ in $V^*$. The subspace
$F$ is mapped by $\sigma$ onto some invariant subspace $U_F\subseteq
V\otimes V^*$. Since $f$ is irreducible, now we can apply the 
theorem~\mythetheorem{5.3}. It yields $U_F=V\otimes\tilde W$, where
$\tilde W\subseteq V^*$ is a subspace in $V^*$.\par
     If we assume that $\tilde W\neq V^*$, then there is some vector
$\bold x\neq 0$ in $V$ such that $\lambda(\bold x)=0$ for all $\lambda
\in V^*$. Applying this fact to $F$ and taking into account 
\mythetag{5.11} and \mythetag{5.12}, we find that this vector $\bold x
\neq 0$ belongs to the kernel of any operator $A\in F$. But the operators
$f(g)\in F$ are non-degenerate, their kernels are zero. This contradiction
shows that $\tilde W=V^*$ and $U_F=V\otimes V^*$. Due to the isomorphism
$\sigma$ then we derive $F=\End(V)$. 
\qed\enddemo
\head
\SectionNum{6}{27} Irreducible representations of the direct product
of groups.
\endhead
\rightheadtext{\S\,6. Representations of the direct product \dots}
     The direct product is the simplest construction for building 
new groups from those already available. \pagebreak Let's recall 
that the group $G_1\times G_2$ is the set of ordered pairs $(g_1,g_2)$ 
with the multiplication rule
$$
\align
(g_1,g_2)\,\cdot\,&(\tilde g_1,\tilde g_2)=(g_1\cdot\tilde g_1,g_2
\cdot\tilde g_2),\\
&\text{where \ }g_1,\tilde g_1\in G_1
\text{\ \ and \ }g_2,\tilde g_2\in G_2.
\endalign
$$
The construction of direct product of groups is in a good agreement
with the construction of tensor product of their representations. Let
$(f_1,G_1,V_1)$ and $(f_2,G_2,V_2)$ are representations of the groups
$G_1$ and $G_2$ respectively. Let's define a representation of the
group $G=G_1\times G_2$ in the space $V_1\otimes V_2$ by the formula
$$
\hskip -2em
\gathered
f(g)(\bold x\otimes\bold y)=
f(g_1,g_2)(\bold x\otimes\bold y)=\\
=(f_1(g_1)\bold x)\otimes(f_2(g_2)\bold y).
\endgathered
\mytag{6.1}
$$
It is easy to verify that the definition \mythetag{6.1} is correct.
\mydefinition{6.1} The representation $(f,G_1\times G_2,V_1\otimes 
V_2)$ given by the formula \mythetag{6.1} is called the {\it tensor
product\/} of the representations $(f_1,G_1,V_1)$ and $(f_2,G_2,V_2)$.
It is denoted $f=f_1\otimes f_2$.
\enddefinition
    Note that the earlier construction of the tensor product given 
by the definition~\mythedefinition{5.1} is embedded into the 
construction~\mythedefinition{6.1}. Indeed, let's consider the 
diagonal in the direct product $G\times G$:
$$
D=\{(g_1,g_2)\in G\times G\!:\ g_1=g_2\}.
$$
It is easy to see that $G\cong G$. The restriction of the representation
\mythetag{6.1} to the diagonal $D$ coincides with the representation
\mythetag{5.1}, where $f=f_1$ and $h=f_2$.
\mytheorem{6.1} The tensor product $(f,G_1\times G_2,V_1\otimes V_2)$ of 
two finite-dimensional representations $(f_1,G_1,V_1)$ and $(f_2,G_2,V_2)$
in complex vector spaces $V_1$ and $V_2$  is irreducible if and only if
both multiplicands $f_1$ and $f_2$ are irreducible.
\endproclaim
\demo{Proof} Let's begin with proving the necessity in the formulated 
proposition. Assume that $(f,G_1\times G_2,V_1\otimes V_2)$ is irreducible.
And assume that the irreducibility condition for $(f_1,G_1,V_1)$ and 
$(f_2,G_2,V_2)$ is broken. For the sake of certainty assume that the second
representation $(f_2,G_2,V_2)$ is reducible. Then $f_2$ has a non-trivial
invariant subspace $\{0\}\subsetneq W_2\subsetneq V_2$. But in this case
$V_1\otimes W_2$ is a non-trivial invariant subspace for $f=f_1\otimes f_2$.
This contradicts the irreducibility of the representation $f$. The necessity
is proved.\par
     Let's prove the sufficiency. Assume that $(f_1,G_1,V_1)$ and 
$(f_2,G_2,V_2)$ are irreducible. In order to prove the irreducibility of
$(f,G_1\times G_2,V_1\otimes V_2)$ we use the irreducibility criterion in
form of the theorem~\mythetheorem{3.1}. Let's choose an arbitrary vector
$\bold u\neq 0$ in $V_1\times V_2$ and consider its orbit. The vector
$\bold u$ can be written as
$$
\hskip -2em
\bold u=\bold x_1\otimes\bold y_1+\ldots+\bold x_k\otimes\bold y_k.
\mytag{6.2}
$$
Without loss of generality we can assume that the vectors $\bold y_1,
\,\ldots,\,\bold y_k$ in \mythetag{6.2} are linearly independent. The
expansion \mythetag{6.2} is not unique. However, if the linearly 
independent vectors $\bold y_1,\,\ldots,\,\bold y_k$ are fixed, then
the corresponding vectors $\bold x_1,\,\ldots,\,\bold x_k$ are 
determined uniquely. Without loss of generality we can assume them
to be nonzero.\par
     Now let's apply the theorem~\mythetheorem{5.4}. Let $A\!:\,V_2\to
V_2$ be a linear operator satisfying the following condition:
$$
\hskip -2em
A\bold y_1=y_1,\ A\bold y_2=0,\ \ldots,\ A\bold y_k=0.
\mytag{6.3}
$$
Since the vectors $\bold y_1,\,\ldots,\,\bold y_k$ in \mythetag{6.3}
are linearly independent, such an operator $A$ does exist. Applying 
the theorem~\mythetheorem{5.4} to the representation $f_2$, we conclude
that the operator $A$ belong to the span of the representation
operators, i\.\,e\.
$$
\pagebreak
A=\sum^q_{i=1}\alpha_i\,f_2(g_i)\text{, \ where \ }g_i\in G_2.
$$
Let's apply the operator $1\otimes A$ to the vector \mythetag{6.2}.
This yields
$$
\hskip -2em
(1\otimes A)\bold u=\sum^k_{i=1}\bold x_i\otimes A\bold y_i=
\bold x_1\otimes\bold y_1.
\mytag{6.4}
$$
On the other hand, for the same quantity we get
$$
\hskip -2em
(1\otimes A)\bold u=\sum^q_{i=1}\alpha_i\,(1\otimes f(g_i))\bold u
=\sum^q_{i=1}\alpha_i\,f(e_1,g_i)\bold u.
\mytag{6.5}
$$
Here $e_1$ is the unit element of the group $G_1$. Comparing 
\mythetag{6.4} and \mythetag{6.5}, we see that the vector 
$\bold x_1\otimes\bold y_1$ belongs to the orbit of the vector 
$\bold u$ from \mythetag{6.3}. Due to the irreducibility of the
representation $f_1$ the orbit of the vector $\bold x_1$ spans 
$V_1$. For the similar reasons the orbit of the vector $\bold y_1$ 
spans $V_2$. These facts mean that any two vectors $\bold x\in V_1$ 
and $\bold y\in V_2$ can be obtained as 
$$
\hskip -2em
\aligned
&\bold x=\sum^r_{i=1}\beta_i\,f(g_i)\bold x_1\text{, \ where \ }
g_i\in G_1;\\
&\bold y=\sum^s_{j=1}\gamma_i\,f(g_j)\bold y_1\text{, \ where \ }
g_i\in G_2.
\endaligned
\mytag{6.6}
$$
From \mythetag{6.6} we immediately derive
$$
\bold x\otimes\bold y=\sum^r_{i=1}\sum^s_{j=1}\beta_i\,
\gamma_i\,f(g_i,g_j)(\bold x_1\otimes\bold y_1),
$$
where $(g_i,g_j)\in G_1\times G_2$. Hence, an arbitrary vector of the
form $\bold x\otimes\bold y$ belongs to the orbit of the vector 
$\bold x_1\otimes\bold y_1$ from \mythetag{6.4}, and this vector in 
turn belongs to the orbit of the vector $\bold u$ from \mythetag{6.2}.
However, we know that the vectors of the form $\bold x\otimes\bold y$ 
spans the whole space $V_1\otimes V_2$. As a result we have proved that
the orbit of an arbitrary vector $\bold u\in V_1\otimes V_2$ spans the
whole space of the representation $f=f_1\otimes f_2$. According to the
theorem~\mythetheorem{3.1}, this representation is irreducible. Thus, 
the theorem~\mythetheorem{6.1} is proved.
\qed\enddemo
\mytheorem{6.2} Any finite-dimensional irreducible representation 
$\varphi$ of the direct product of two groups $G_1$ and $G_2$ in 
a complex space $U$ is isomorphic to the tensor product of two 
irreducible representations $(f_1,G_1,V_1)$ and $(f_2,G_2,V_2)$ 
of the groups $G_1$ and $G_2$.
\endproclaim
    Let $\varphi(g_1,g_2)$ be the representation operators for the 
representation $\varphi$ of the group $G_1\times G_2$ in the space $U$.
Then the operators of the form $\varphi(g_1,e_2)$, where $e_2$ is the
unit element of the group $G_2$, define a representation of the group
$G_1$. In general case it is reducible. Let $V_1\subseteq U$ be some 
irreducible subspace in $U$. Denote
$$
\hskip -2em
\varphi_1(g_1)=\varphi(g_1,e_2)\text{, \ where \ }g_1\in G_1.
\mytag{6.7}
$$
The restrictions of the operators \mythetag{6.7} to $V_1$ define some
irreducible representation of the group $G_1$. We denote them as
$$
f_1(g_1)=\varphi_1(g_1)\,\hbox{\vrule height 8pt depth 10pt
width 0.5pt}_{\,V_1}=\varphi_1(g_1,e_2)\,\hbox{\vrule height 8pt depth 10pt
width 0.5pt}_{\,V_1}\text{, \  where \ } g_1\in G_1.  
$$
By analogy to \mythetag{6.7} we introduce the following operators defining
some representation of the group $G_2$ in $U$:
$$
\hskip -2em
\varphi_2(g_2)=\varphi(e_1,g_2)\text{, \ where \ }g_2\in G_2.
\mytag{6.8}
$$
Then we denote by $F_2$ the span of the set of all operators \mythetag{6.8}.
It is a subspace in the space of the operators $\End(U)$:
$$
F_2=\left<\{\varphi_2(g_2)\!:\ g_2\in G_2\}\right>.
$$
The operators from $F_2$ commute with all operators \mythetag{6.7} since
the operators $\varphi_2(g_2)$ spanning $F_2$ commute with $\varphi_1(g_1)$.
For each operator $A\in F_2$ we denote by $\tilde A$ the restriction of $A$
to the subspace $V_1$. The operators $\tilde A$ should be treated as the
elements of the linear space $\tilde F_2\subseteq\Hom(V_1,U)$:
$$
\tilde A\!:\,V_1\to U.
$$\par
     The operators $A\in F_2$ and $\tilde A\in\tilde F_2$ deserve a special
consideration. Let's define a subspace $V_A=AV_1=\tilde AV_1=\Img\tilde A
\subseteq U$. Since $A$ commute with operators \mythetag{6.7}, the subspace
$V_A$ is invariant with respect to the operators $\varphi_1(g_1)$.
Therefore we have a representation of the group $G_1$ in $V_A$:
$$
f_A(g_1)=\varphi_1(g_1)\,\hbox{\vrule height 8pt depth 10pt
width 0.5pt}_{\,V_A}=\varphi(g_1,e_2)\,\hbox{\vrule height 8pt 
depth 10pt width 0.5pt}_{\,V_A}\text{, \ where \ }g_1\in G_1.
$$
The mapping $\tilde A\!:\,V_1\to V_A$ interlaces the representations $f_1$
and $f_A$ in $V_1$ and $V_A$. Indeed, we have  
$$
\tilde A\compos f_1(g_1)=A\compos\varphi(g_1)\,\hbox{\vrule height 8pt
depth 10pt width 0.5pt}_{\,V_1}=\varphi(g_1)\compos A\,\hbox{\vrule 
height 8pt depth 10pt width 0.5pt}_{\,V_1}=f_A(g_1)\compos\tilde A.
$$
The mapping $\tilde A\!:\,V_1\to V_A$ is surjective by its definition. 
The kernel $\Ker\tilde A\subseteq V_1$ of this mapping is invariant with
respect to the representation $f_1$. Since $f_1$ is irreducible, we have
two mutually excluding options:
$$
\aligned
\Ker\tilde A=V_1\ &\Rightarrow\ \tilde A=0
\text{\ \ and \ }V_A=\{0\};\\
\Ker\tilde A=\{0\}\ &\Rightarrow\ \tilde A\text{\ \ is an
isomorphism and \ }f_1\cong f_A.
\endaligned\quad
\mytag{6.9}
$$
Let's study the second option in \mythetag{6.9}. Denote $W=V_1\cap V_A$.
The operators $f_1(g_1)$ and $f_A(g_1)$ upon restricting to $W$ do coincide.
Therefore $W\subseteq V_1$ is invariant with respect to $f_1$. Applying the
irreducibility of $f_1$ again, we get the following two options:
$$
\aligned
W=V_1\ &\Rightarrow\ V_A=V_1\text{\ \ and \ }f_A=f_1;\\
W=\{0\}\ &\Rightarrow\ V_A\cap V_1=\{0\}\text{\ \ and \ }f_A\cong f_1.
\endaligned\quad
\mytag{6.10}
$$
Combining \mythetag{6.9} and \mythetag{6.10}, we find
$$
\align
&V_A=\{0\},\quad f_A=0,\quad\tilde A=0;\\
&V_A=V_1,\quad f_A=f_1,\quad\tilde A=\lambda\cdot 1;
\mytag{6.11}\\
&V_A\cap V_1=\{0\},\quad f_A\cong f_1,\quad\tilde A\text{\ \ is an
isomorphism.}\qquad
\endalign
$$
The condition $\tilde A=\lambda\cdot 1$ in the second option of 
\mythetag{6.11} follows from Schur's lemma~\mythelemma{5.2}.\par
    Let $\bold u$ be some nonzero vector in $V_1$. We fix this vector 
and consider the subspace $V_2\subseteq U$ obtained by applying the
operators $A\in F_2$ upon the vector $\bold u$:
$$
V_2=F_2\bold u=\{\bold v\in U\!:\ \bold v=A\bold u\text{\ \ for some \ }
A\in F_2\}.\quad
\mytag{6.12}
$$
The subspace $V_2$ is invariant with respect to the operators
\mythetag{6.8}. Therefore we have a representation $(f_2,G_2,V_2)$
of the group $G_2$. It is given by the operators
$$
\hskip -2em
f_2(g_2)=\varphi_2(g_2)\,\hbox{\vrule height 8pt depth 10pt
width 0.5pt}_{\,V_2}=\varphi_1(e_1,g_2)\,\hbox{\vrule height 8pt depth 10pt
width 0.5pt}_{\,V_2}.  
\mytag{6.13}
$$\par
     Due to the definition \mythetag{6.12} for any vector $\bold y
\in V_2$ there is a mapping $\tilde A\in\tilde F_2$ such that 
$\bold y=\tilde A\bold u$. Let's prove that such a mapping is uniquely
fixed by the vector $\bold y\in V_2$. According to \mythetag{6.11}, we
study three possible options.\par
     If $\bold y=0$, then $\Ker\tilde A\neq 0$. Due to \mythetag{6.9}
the only operator $\tilde A\in F_2$ satisfying the condition 
$\bold y=\tilde A\bold u$ is the identically zero mapping $\tilde A=0$.
This case corresponds to the first option in \mythetag{6.11}.\par
     If $\bold y\neq 0$ and $\bold y\in V_1$, then from $\bold y=
\tilde A\bold u$ we derive that $\bold y\in V_1\cap V_A$. Hence the
intersection $V_1\cap V_A$ is nonzero and we have the first option 
in \mythetag{6.10}, which is equivalent to the second option of
\mythetag{6.11}. Hence, $\tilde A=\lambda\cdot 1$ and $\bold y=\lambda
\bold u$. The number $\lambda$ relating two collinear vectors is 
uniquely fixed by these two vectors. Therefore, the mapping 
$\tilde A=\lambda\cdot 1$ is also unique.\par
     And finally, the third case, where $\bold y\notin V_1$. Due to
\mythetag{6.11} in this case we have $V_1\cap V_A=\{0\}$ and the 
mapping $\tilde A\!:\,V_1\to V_A$ is bijective. Assume for a while 
that the condition $\bold y=\tilde A\bold u$ does not fix the mapping
$\tilde A\in F_2$ uniquely. Let $\tilde A_1$ and $\tilde A_2$ be two
such mappings. Their associated subspaces $V_{A_1}$ and $V_{A_2}$ do
coincide. Indeed, $\bold y\in V_{A_1}\cap V_{A_2}\neq\{0\}$. Hence,
$V_{A_1}\cap V_{A_2}$ is a non-trivial invariant subspace for the
irreducible representations $f_{A_1}\cong f_1$ and $f_{A_2}\cong f_1$. 
So, $V_{A_1}\cap V_{A_2}=V_{A_1}=V_{A_2}$. Using $V_{A_1}=V_{A_2}$ and
the bijectivity of the mappings 
$$
\xalignat 2
&\tilde A_1\!:V_1\to V_{A_1},
&&\tilde A_2\!:V_1\to V_{A_2},
\endxalignat
$$
we invert one of them and consider the operator 
$\tilde A_3=\tilde A_2^{-1}\compos\tilde A_1$. This is a non-degenerate 
operator in $V_1$. It implements the automorphism of the representation 
$f_1$, i\.\,e\. it interlaces the operators $f_1(g_1)$ with themselves:
$$
\tilde A_3\compos f_1(g_1)=f_1(g_1)\compos\tilde A_3.
$$
Using the irreducibility of $f_1$ and applying Schur's 
lemma~\mythelemma{5.2}, we get $\tilde A_3=\lambda\cdot 1$. This 
yields $\tilde A_1=\lambda\tilde A_2$. Now from the conditions 
$\bold y=\tilde A_1\bold u$ and $\bold y=\tilde A_2\bold u$ we
derive $\lambda=1$. Hence, $\tilde A_2=\tilde A_1$. Thus, the
uniqueness of $\tilde A$ is established.\par
     For the mapping $\tilde A\in\tilde F_2$, which we uniquely 
determine from the condition $\bold y=\tilde A\bold u$, we use
the notation $\tilde A=\tilde A(\bold y)$. The dependence of 
$\tilde A$ on the vector $\bold y$ can be treated as a mapping
$\tilde A\!:\,V_2\to\Hom(V_1,U)$. It is easy to verify that
this mapping is linear. It satisfies the equality
$$
\tilde A(f_2(g_2)\bold y)=\varphi_2(g_2)\compos\tilde A(\bold y),
\mytag{6.14}
$$
where the operator $f_2(g_2)$ is given by \mythetag{6.13}. Let's prove 
the equality \mythetag{6.14}. Remember that $\tilde A(\bold y)$ is
the restriction to $V_1$ of some operator $A_1\in F_2$ such that
$$
A_1\bold u=\tilde A(\bold y)\bold u=\bold y.
$$
But the operator $A_2=\varphi_2\compos A_1$ also belongs to $F_2$ 
(see the definition of the space $F_2$ above). For $A_2$ we derive
$$
\pagebreak
A_2\bold u=\varphi_2(g_2)A_1\bold u=\varphi_2(g_2)\bold y=
f_2(g_2)\bold y.
$$
Therefore the restriction of $A_2$ to $V_1$ coincides with 
$\tilde A(f_2(g_2)\bold y)$. The equality \mythetag{6.14} is
proved.\par
     The next step in proving the theorem~\mythetheorem{6.2} is to
apply the mapping $\tilde A(\bold y)$ for building the isomorphism
of the representation $f=f_1\otimes f_2$ and the representation 
$\varphi$. But before doing it note that we have no information 
on whether the representation
$f_2$ in \mythetag{6.13} is irreducible or not. Fortunately we can 
assume $f_2$ to be irreducible due to the following reasons. Let 
$\tilde V_2\subseteq V_2$ be some irreducible invariant subspace 
for the representation \mythetag{6.13}. If $\bold u\in\tilde V_2$, then
$\tilde V_2=V_2$. This fact follows from the theorem~\mythetheorem{3.1}.
In the case, where $\bold u\notin\tilde V_2$, we choose a nonzero vector
$\tilde\bold u$ and take a mapping $\tilde A\in\tilde F_2$ such that
$\tilde A\bold u=\tilde\bold u$. We have already proved the existence and
uniqueness of such a mapping $\tilde A=\tilde A(\tilde\bold u)$. In our 
case $\tilde A\neq 0$. Therefore, due to \mythetag{6.11} we see that the 
mapping $\tilde A$ is bijective, it establishes the isomorphism of 
representations $f_1\cong f_A$. Because of the isomorphism $f_1\cong f_A$
we can replace $f_1$ by $f_A$, which is also irreducible. The latter
representation is preferable since its space $V_A$ comprises the vector
$\tilde\bold u$. The orbit of the vector $\tilde\bold u$ spans the 
irreducible subspace $\tilde V_2$ within the space of the representation 
$\varphi_2$. For this reason we should come back to the beginning of
our constructions and assume that $V_1$ is exactly that irreducible
subspace of $\varphi_1$ which comprises some vector $\bold u$ generating
an irreducible subspace of the representation $\varphi_2$. Just above
we have demonstrated that such a choice of the subspace $V_1$ is
possible.\par
    Thus, under a proper choice of the subspace $V_1$ both representations 
$(f_1,G_1,V_1)$ and $(f_2,G_2,V_2)$ are irreducible. We consider their
tensor product $f=f_1\otimes f_2$ and then construct the mapping 
$\sigma\!:\,V_1\otimes V_2\to U$ by means of the following formula:
$$
\hskip -2em
\sigma(\bold x\otimes\bold y)=\tilde A(\bold y)\bold x
\text{, \ where \ }\bold x\in V_1,\ \bold y\in V_2.
\mytag{6.15}
$$
Let's show that the mapping \mythetag{6.15} is an interlacing mapping for
the representations $(f_1\otimes f_2,G_1\times G_2,V_1\otimes V_2)$ and
$(\varphi,G_1\times G_2,U)$. Indeed, we easily derive
$$
\aligned
\hskip -2em
&\gathered
\varphi(g_1,g_2)\,\sigma(\bold x\otimes\bold y)=\varphi_1(g_1)
\varphi_2(g_2)\tilde A(\bold y)\bold x=\\
=\varphi_2(g_2)\tilde A(\bold y)\,\varphi_1(g_1)\bold x,
\endgathered\\
\vspace{4ex}
&\gathered
\sigma f(g_1,g_2)(\bold x\otimes\bold y)=\sigma((f_1(g_1)\bold x)
\otimes(f_2(g_2)\bold y))=\\
=\tilde A(f_2(g_2)\bold y)f_1(g_1)\bold x.
\endgathered
\endaligned
\mytag{6.16}
$$
The values of the right hand sides in two above formulas \mythetag{6.16} 
do coincide due to \mythetag{6.14}. Therefore, from \mythetag{6.16} we 
extract
$$
\varphi(g_1,g_2)\compos\sigma=\sigma\compos f(g_1,g_2).
$$
This is exactly the equality \mythetag{1.1} written for the representations
$f$ and $\varphi$. The mapping $\sigma$ implements an isomorphism of these 
two representations. Note that 
$$
\sigma(\bold u\times\bold u)=\tilde A(\bold u)\bold u=\bold u\neq 0.
$$
Therefore $\sigma\neq 0$. Now it is sufficient to use the irreducibility 
of representations $f=f_1\otimes f_2$ and $\varphi$. Applying Schur's 
lemma~\mythelemma{5.1}, we conclude that $\sigma$ is an isomorphism. 
The irreducibility of $f$ is derived from the irreducibility of $f_1$ 
and  $f_2$ due to the previous theorem. Thus, the proof of the 
theorem~\mythetheorem{6.2} is completed.
\head
\SectionNum{7}{36} Unitary representations.
\endhead
\rightheadtext{\S\,7. Unitary representations.}
\mydefinition{7.1} A finite-dimensional complex linear vector space $V$
equipped with a symmetric positive sesquilinear form is called a {\it
Hermitian space}.
\enddefinition
     Let's recall that a {\it sesquilinear form\/} in $V$ is a
complex-valued numeric function $\varphi(\bold x,\bold y)$ with two
vectorial arguments $\bold x,\bold y\in V$ such that it satisfies
the following four conditions:
\roster
\item\quad $\varphi(\bold x_1+\bold x_2,\bold y)=\varphi(\bold x_1,\bold y)
       +\varphi(\bold x_2,\bold y)$;
\item\quad $\varphi(\alpha\,\bold x,\bold y)=\overline\alpha\,
       \varphi(\bold x,\bold y)$;
\item\quad $\varphi(\bold x,\bold y_1+\bold y_2)\varphi(\bold x,\bold y_1)
       +\varphi(\bold x,\bold y_2)$;
\item\quad $\varphi(\bold x,\alpha\,\bold y)=\alpha\,
       \varphi(\bold x,\bold y)$.
\endroster
     The bar sign over $\alpha$ in the second condition is the complex
conjugation sign. The conditions \therosteritem{1}--\therosteritem{4}
are usually are usually complemented with the conditions of symmetry 
and positivity:
\roster
\item[5]\quad $\varphi(\bold x,\bold y)
       =\overline{\varphi(\bold y,\bold x)}$;
\item\quad $\varphi(\bold x,\bold x)>0\text{\ \ for all \ }
       \bold x\neq 0$.
\endroster
     The condition \therosteritem{5} implies that $\varphi(\bold x,
\bold x)$ is a real number. The condition \therosteritem{6} strengthens
condition \therosteritem{5} requiring $\varphi(\bold x,\bold x)$ to be
a positive number. A form $\varphi(\bold x,\bold y)$ is called 
{\it non-degenerate} if $\varphi(\bold x,\bold y)=0$ for all $\bold y\in V$
implies $\bold x=0$. Note that the positivity of a form implies its
non-degeneracy.\par
     The symmetric positive form declared in the definition of a Hermitian
space is called the {\it Hermitian scalar product}. For this form we fix
the following notation:
$$
\varphi(\bold x,\bold y)=\langle\bold x|\bold y\rangle.
$$\par
     Let $\bold e_1,\,\ldots,\,\bold e_n$ be a basis in a space $V$. The
quantities $g_{ij}=\langle\bold e_i|\bold e_j\rangle$ compose the Gram
matrix of the basis $\bold e_1,\,\ldots,\,\bold e_n$. They satisfy the
relationship $g_{ij}=\overline{g_{ji}\vphantom{\vrule height 5pt}}$.
It follows from the symmetry of the scalar product.\par
     A basis, the Gram matrix of which is the unit matrix, is called an
{\it orthonormal basis}. Orthonormal bases do exist because each symmetric
sesquilinear in a finite-dimensional space is diagonalizable.
\mydefinition{7.2} A linear operator $A\!:\,V\to V$ in a Hermitian space 
$V$ is called a {\it Hermitian operator\/} if $\langle\bold x|A\bold y
\rangle=\langle A\bold x|\bold y\rangle$ for any two vectors $\bold x$ 
and $\bold y$ in $V$.
\enddefinition
     There is the standard theory of Hermitian operators in
finite-dimensional Hermitian spaces. We give basic facts of this theory
without proofs for the reader to recall them.
\mytheorem{7.1} Hermitian operators of a finite-dimensional\linebreak 
Hermitian space are in a one-to-one correspondence with symmetric
sesquilinear forms:
$$
\hskip -2em
\varphi_A(\bold x,\bold y)=\langle\bold x|A\bold y\rangle.
\mytag{7.1}
$$
Non-degenerate operators correspond to non-degenerate forms.
\endproclaim
\mydefinition{7.3} A Hermitian operator $A$ is called a {\it positive 
operator\/} if the corresponding form $\varphi_A$ is positive.
\enddefinition
\mytheorem{7.2} Each Hermitian operator $A$ is diagonalizable, its 
eigenvalues are real numbers, and eigenvectors corresponding to 
distinct eigenvalues $\lambda_i\neq\lambda_j$ are perpendicular to 
each other.
\endproclaim
\mytheorem{7.3} An operator $A$ is a Hermitian operator if and only 
if its eigenvalues are real numbers and if is diagonalizes in some
orthonormal basis.
\endproclaim
    The proofs of the theorems~\mythetheorem{7.1}, \mythetheorem{7.2}, 
and \mythetheorem{7.3} can be found in many standard textbooks on
linear algebra. Apart from them, we need one more theorem, which also 
can be found in some textbooks, but it is less standard.
\mytheorem{7.4} Let $A$ be a diagonalizable operator such that its
eigenvalues $\lambda_1,\,\ldots,\,\lambda_n$ are real non-negative 
numbers. Then there is a unique operator $B$ with eigenvalues $\mu_i
\geqslant 0$ such that $B^2=A$ and $B$ commutes with any operator $C$
that commutes with $A$. If $A$ is a Hermitian operator, then the
corresponding operator $B$ is a Hermitian operator too.
\endproclaim
    The operator $B$ declared in the theorem~\mythetheorem{7.4} is
naturally called the {\it square root\/} of the operator $A$. Let's
prove its existence. Let $\bold e_1,\,\ldots,\,\bold e_n$ be a basis
composed by eigenvectors of the operator $A$ corresponding to its
eigenvalues $\lambda_1,\,\ldots,\,\lambda_n$. The operator $B$ is
defined through its action upon basis vectors:
$$
\hskip -2em
B\bold e_i=\sqrt{\lambda_i}\,\bold e_i,\qquad i=1,\,\ldots,\,n.
\mytag{7.2}
$$
Due to this definition the operator $B$ is diagonalized in the same 
basis as the operator $A$, its eigenvalues $\mu_i=\sqrt{\lambda_i}$
are real and non-negative numbers.\par
     Let's study the problem of commuting for the operators $B$ and $C$.
If the operator $A$ commutes with $C$, this means that 
$$
\hskip -2em
(\lambda_i-\lambda_j)\,C^i_j=0,
\mytag{7.3}
$$
where $C^i_j$ are the matrix elements for the operator $C$ in the 
basis $\bold e_1,\,\ldots,\,\bold e_n$. The condition \mythetag{7.3} 
is equivalent to $C^i_j=0$ for all $\lambda_i\neq\lambda_j$. But
$\lambda_i\neq\lambda_j$ implies $\mu_i\neq\mu_j$. Therefore the
operator $B$ defined by the formula \mythetag{7.2} commutes with
any operators $C$ that commutes with $A$. If the operator $A$ is a
Hermitian operator, then the basis $\bold e_1,\,\ldots,\,\bold e_n$
can be chosen to be an orthonormal basis. In this case, applying 
the theorem~\mythetheorem{7.3}, we find that $B$ is a Hermitian
operator too.\par
    Now we need to prove the uniqueness of the operator $B$ declared
in the theorem~\mythetheorem{7.4}. The condition $C^i_j=0$ for all
$\lambda_i\neq\lambda_j$, which follows from $A\compos C=C\compos A$,
can be formulated in an invariant (basis-free) way.
\myproposition{7.1} An operator $C$ commutes with a diagonalizable
operator $A$ if and only if all eigenspaces of the operator $A$ are
invariant with respect to the operator $C$.
\endproclaim
     Under the assumptions of the theorem~\mythetheorem{7.4} let's
take $C=A$ and apply the proposition~\mytheproposition{7.1} to the
operator $B$. From $B\compos A=A\compos B$ in this case we derive that
all eigenspaces of the operator $A$ are invariant under the action of 
the operator $B$. The requirement that $B$ is diagonalizable now means
that both $A$ and $B$ can be diagonalized simultaneously in some basis.
The conditions $B^2=A$ and $\mu_i\geqslant 0$ then fix the unique choice 
of the operator $B$ defined by the relationships \mythetag{7.2}.\par
\mydefinition{7.4} A linear mapping $T\!:\,V\to W$ from some Hermitian
space $V$ to another Hermitian space $W$ is called an {\it isometry\/} 
if $\langle T\bold x|T\bold y\rangle=\langle\bold x|\bold y\rangle$ for 
all $\bold x,\bold y\in V$, i\.\,e\. if it preserves the scalar product.
\enddefinition
    Due to the non-degeneracy of the sesquilinear forms determining 
scalar products in $V$ and $W$ each isometry $T\!:\,V\to W$ is an 
injective mapping.
\mydefinition{7.5} A linear operator $T\in\End(V)$ is called a {\it 
unitary operator\/} if it implements an isometry $T\!:\,V\to V$.
\enddefinition
    Unitary operators are non-degenerate. Their determinants and their 
eigenvalues satisfy the following relationships:
$$
\xalignat 2
&|\det T|=1,&&|\lambda|=1.
\endxalignat
$$
Unitary operators in a Hermitian space $V$ form the group $\MatGrU(V)$,
which is a subgroup in $\Aut(U)$. Unitary operators with unit determinant,
in turn, form the group $\MatGrSU(V)\subsetneq\MatGrU(V)$.
\mydefinition{7.6} A representation $(f,G,V)$ of a group $G$ in 
a Hermitian space $V$ is called a {\it unitary representation\/} 
if all operators of this representation $f(g)$ are unitary operators.
\enddefinition
     Unitary representations constitute an important subclass in the class 
of general representations of groups. First of all this is because unitary
representations arise in applications of the theory of representations to 
the quantum mechanics. The role of the following useful fact is also
substantial.
\mytheorem{7.5} Each unitary representation $(f,G,V)$ is completely
reducible.
\endproclaim 
     Indeed, let $U\subseteq V$ be an invariant subspace for the operators
of the representation $f$. In the case of a unitary operator $f(g)$ the
orthogonal complement to an invariant subspace
$$
U_{\sssize\perp}=\{\bold x\in V\!:\ \langle\bold x|\bold y\rangle=0
\text{\ \ for all \ }\bold y\in U\}
$$
is also an invariant subspace. The subspaces $U$ and $U_{\sssize\perp}$
intersect trivially (i\.\,e\. at zero vector only), their direct sum
coincides with $V$. Therefore, $U_{\sssize\perp}$ is a required invariant
direct complement for $U$. The complete reducibility of $f$ is shown.
\mycorollary{7.1} Each representation $f$ which is equivalent to some
unitary representation $h$ is completely reducible.
\endproclaim 
     Let $A\!:\,V\to W$ be an interlacing mapping which implements the
isomorphism of $f$ and $h$. Then each invariant subspace $U$ of $f$ has
the invariant direct complement $A^{-1}((AU)_{\sssize\perp})$.\par
     Along with the concept of equivalence, in the class of unitary
representations we have the concept of {\it unitary equivalence}.
\mydefinition{7.7} Two unitary representations $(f,G,V)$ and $(h,G,W)$
are called unitary equivalent, if there is an isometry $A\!:\,V\to W$ 
implementing an isomorphism of them.
\enddefinition
    The following theorem shows that despite the difference in definitions
the concepts of equivalence and unitary equivalence do coincide.
\mytheorem{7.6} If two unitary representations $f$ and $h$ are equivalent, 
then they are unitary equivalent.
\endproclaim 
    In order to prove this theorem we need an auxiliary fact which is
formulated as a lemma.
\mylemma{7.1} Let $A\!:\,V\to W$ be a bijective linear mapping from a
Hermitian space $V$ to another Hermitian space $W$. Then it can be expanded
as a composition $A=T\compos B$, where $T\!:\,V\to W$ is an isometry and
$B$ is a positive Hermitian operator in $V$.
\endproclaim
\demo{Proof of the lemma} Let's consider the following sesquilinear form 
in the space $V$:
$$
\hskip -2em
\varphi(\bold x,\bold y)=\langle A\bold x|A\bold y\rangle.
\mytag{7.4}
$$
It is easy to see that the form \mythetag{7.4} is symmetric and positive.
We apply the theorem~\mythetheorem{7.1} in order to define a Hermitian 
positive operator $D$ in $V$. The associated sesquilinear form \mythetag{7.1} 
of this operator coincides with \mythetag{7.4}. This condition yields
$$
\hskip -2em
\langle\bold x|D\bold y\rangle=\langle A\bold x|A\bold y\rangle.
\mytag{7.5}
$$
Using the operator $D$ and applying the theorem~\mythetheorem{7.4} to it,
we construct a positive Hermitian operator $B$ being the square root of $D$,
i\.\,e\. $B^2=D$. Now let's consider a mapping $T\!:\,V\to W$ defined as
the composition $T=A\compos B^{-1}$. Note that $B$ is non-degenerate since
it is positive. Therefore it is invertible. The rest is to show that 
$T$ is an isometry. Indeed, we have
$$
\langle T\bold x|T\bold y\rangle=\langle A\compos B^{-1}\bold x|
A\compos B^{-1}\bold y\rangle=\langle B^{-1}\bold x|D\compos B^{-1}
\bold y\rangle.\quad
\mytag{7.6}
$$
The last equality in the chain \mythetag{7.6} is provided by 
\mythetag{7.5}. The further calculations are obvious:
$$
\langle B^{-1}\bold x|D\compos B^{-1}\bold y\rangle
=\langle B^{-1}\bold x|\compos B\bold y\rangle
=\langle B\compos B^{-1}\bold x|\bold y\rangle
=\langle\bold x|\bold y\rangle.
$$
Combining this equality with \mythetag{7.6}, we get
$\langle T\bold x|T\bold y\rangle=\langle\bold x|\bold y\rangle$ 
for any two vectors $\bold x,\bold y\in V$. Hence, $T$ is an 
isometry. The lemma is proved.
\qed\enddemo
\demo{Proof of the theorem~\mythetheorem{7.6}} The mapping $A=T\compos B$
in this case implements an isomorphism of two unitary representations
$f$ and $h$. Therefore, we have 
$$
\hskip -2em
T\compos B\compos f(g)=h(g)\compos T\compos B
\text{\ \  for all \ }g\in G.
\mytag{7.7}
$$
Let's show that the operator $B$ commutes with $f(g)$. For this purpose 
we show that $D=B^2$ commutes with $f(g)$:
$$
\langle\bold x|f(g)B^2\bold y\rangle=\langle f(g^{-1})\bold x|B
B\bold y\rangle=\langle Bf(g^{-1})\bold x|B\bold y\rangle.
$$
Here we used the facts that $f(g)$ is a unitary operator and $B$ is a
Hermitian operator. Let's continue our calculations using the isometry
of the mapping $T$:
$$
\langle Bf(g^{-1})\bold x|B\bold y\rangle=\langle Bf(g^{-1})\bold x|
TB\bold y\rangle=\langle TBf(g^{-1})\bold x|B\bold y\rangle.
$$
But $T\compos B\compos f(g^{-1})=h(g^{-1})\compos T\compos B$. This fact
follows from \mythetag{7.7}. Taking into account this equality and taking
into account that $h(g)$ is a unitary operator, we get
$$
\langle TBf(g^{-1})\bold x|B\bold y\rangle=\langle h(g^{-1})TB\bold x
|B\bold y\rangle=\langle TB\bold x|h(g)TB\bold y\rangle.
$$
Now let's use again the relationship \mythetag{7.7} written as
$h(g)\compos T\compos B=T\compos B\compos h(g)$. Then we take 
into account the isometry of $T$:
$$
\langle TB\bold x|h(g)TB\bold y\rangle=\langle TB\bold x|TBf(g)
\bold y\rangle=\langle T^{-1}TB\bold x|Bf(g)\bold y\rangle.
$$
In order to complete this series of calculations, we remember that
$B$ is a Hermitian operator:
$$
\gathered
\langle T^{-1}TB\bold x|Bf(g)\bold y\rangle=\langle B\bold x|Bf(g)
\bold y\rangle=\\
=\langle\bold x|BBf(g)\bold y\rangle=
\langle\bold x|B^2f(g)\bold y\rangle.
\endgathered
$$
As a result we have got $\langle\bold x|f(g)B^2\bold y\rangle
=\langle\bold x|B^2f(g)\bold y\rangle$. Since $\bold x$ and $\bold y$
are arbitrary two vectors, we conclude that the operators $f(g)$ commute 
with the operator $D=B^2$. But the positive Hermitian operator $B$ is
a square root of the positive Hermitian operator $D$. Therefore $B$
commutes with all operators that commute with the operator $D$ (see 
theorem~\mythetheorem{7.4}). As a result we get $f(g)\compos B=B\compos 
f(g)$. Substituting this equality into \mythetag{7.7}, we derive
$T\compos f(g)\compos B=h(g)\compos T\compos B$. Canceling the
non-degenerate operator $B$, we find
$$
\hskip -2em
T\compos f(g)=h(g)\compos T\text{\ \ for all \ }g\in G.
\mytag{7.8}
$$
The equality \mythetag{7.8} means that the isometric mapping $T$
implements an isomorphism of the unitary representations $f$ 
and $h$. Hence, the representations $f$ and $h$ are unitary equivalent.
Thus, the theorem~\mythetheorem{7.6} is proved.
\qed\enddemo
\newpage
\setfirstpage
\topmatter
\title\chapter{2}
Representations of finite groups
\endtitle
\endtopmatter
\leftheadtext{CHAPTER \uppercase\expandafter{\romannumeral 2}.
REPRESENTATIONS OF FINITE GROUPS.}
\document
\chapternum=2
\head
\SectionNum{1}{44} Regular representations of finite groups.
\endhead
\rightheadtext{\S\,1. Regular representations of finite groups.}
     Let $G$ be a finite group and $N=|G|$ be the number of elements 
in this group. Let's consider the set of complex-valued numeric functions on
$G$. We denote it by $L_2(G)$. It is clear that $L_2(G)$ is a complex 
linear vector space of the dimension $\dim(L_2(G))=N$. Let's equip $L_2(G)$
with the structure of a Hermitian space. For this purpose we consider the
scalar product of two functions $u(g)$ and $v(g)$ given by the formula
$$
\hskip -2em
\langle u|v\rangle=\frac{1}{N}\sum_{g\in G}\overline{u(g)}\,v(g).
\mytag{1.1}
$$
Now we define an action of the group $G$ in $L_2(G)$ by defining linear
operators $R(g)\!:\,L_2(G)\to L_2(G)$. Let's set
$$
R(g)v(a)=v(a\,g)\text{\ \ for all \ } a,g\in G\text{\ \ and \ }
v\in L_2(G).\quad
\mytag{1.2}
$$
The operators $R(g)$ act upon functions of $L_2(g)$ by means of right
shifts of their arguments. It is easy to verify that these operators
satisfy the following relationship:
$$
R(g_1)\compos R(g_2)=R(g_1\,g_2).
$$
Hence, the operators $R(g)$, which  act according to \mythetag{1.2},
form a representation of the group $G$ in the space $L_2(G)$. This 
representation is called the {\it right regular representation\/} 
of the group $G$.\par
    Along with the right regular representation there is the {\it left 
regular representation\/} $(L,G,L_2(G))$ of the group $G$. Its operators
are defined as follows:
$$
L(g)v(a)=v(g^{-1}\,a)\text{\ \ for all \ } a,g\in G\text{\ \ and \ }
v\in L_2(G).\quad
\mytag{1.3}
$$
\mytheorem{1.1} The right regular representation defined by the formula
\mythetag{1.2} and the left regular representation defined by the formula
\mythetag{1.3} are unitary representations with respect to the Hermitian
structure given by the scalar product \mythetag{1.1}.
\endproclaim
\demo{Proof} Let's verify by means of direct calculations that $R(g)$ and
$L(g)$ are unitary operators. Assume that $u$ and $v$ are two arbitrary
functions from $L_2(G)$. Then we have
$$
\align
&\gathered
\langle R(g)u|R(g)v\rangle=\frac{1}{N}\sum_{a\in G}
\overline{u(a\,g)}\,v(a\,g)=\\
=\frac{1}{N}\sum_{b\in G}
\overline{u(b)}\,v(b)=\langle u|v\rangle;
\endgathered\\
\vspace{2ex}
&\gathered
\langle L(g)u|L(g)v\rangle=\frac{1}{N}\sum_{a\in G}
\overline{u(g^{-1}\,a)}\,v(g^{-1}\,a)=\\
=\frac{1}{N}\sum_{b\in G}
\overline{u(b)}\,v(b)=\langle u|v\rangle;
\endgathered
\endalign
$$
Here we used the fact that the right shift $a\mapsto b=a\,g$ and the
left shift $a\mapsto g^{-1}\,a$ are two bijective mappings of the group
$G$ onto itself.
\qed\enddemo
\mytheorem{1.2} The right regular representation and the left regular 
representation are unitary equivalent to each other.
\endproclaim
\demo{Proof} In order to prove the theorem it is necessary to 
construct the unitary operator $A\!:\,L_2(G)\to L_2(G)$ interlacing 
the representations $(R,G,L_2(G))$ and $(L,G,L_2(G))$. We define 
this operator as follows:
$$
\hskip -2em
Av(g)=v(g^{-1})\text{\ \ for all \ }g\in G\text{\ \ and \ }
v\in L_2(G).
\mytag{1.4}
$$
The fact that $A$ is a unitary operator with respect to the Hermitian
structure \mythetag{1.1} is shown by the following calculations:
$$
\langle Au|Av\rangle=\frac{1}{N}\sum_{a\in G}
\overline{u(a^{-1})}\,v(a^{-1})=\frac{1}{N}
\sum_{b\in G}\overline{u(b)}\,v(b)=\langle u|v\rangle.
$$
In carrying out these calculations we used the fact that the inversion
operation $a\mapsto b=a^{-1}$ is a bijective mapping of the group $G$
onto itself.\par
     Now let's show that the operator $A$ introduced by means of the
formula \mythetag{1.4} interlaces right and left regular representations.
Let $v$ be some arbitrary function from $L_2(G)$. Assume that $u=L(g)v$
and $w=Av$. Then we have
$$
\gather
AL(g)v(a)=Au(a)=u(a^{-1})=v(g^{-1}\,a^{-1})=\\
=v((a\,g)^{-1})=w(a\,g)=R(g)w(a)=R(g)Av(a).
\endgather
$$
Since $v$ is an arbitrary function from $L_2(G)$, this sequence of
calculations shows that $A\compos L(g)=R(g)\compos A$. The proof 
is over.
\qed\enddemo
\head
\SectionNum{2}{46} Invariant averaging over a finite group.
\endhead
\rightheadtext{\S\,2. Invariant averaging \dots}
     In previous section we have noted that the idea of considering 
numeric functions on a finite group is fruitful enough. The finiteness
of a group $G$ provides the opportunity to define the operation of 
invariant averaging for such functions. For an arbitrary function
$v\in L_2(G)$ we denote by $M[v]$ the number determined by the following
relationship:
$$
\hskip -2em
M[v]=\frac{1}{N}\sum_{g\in G}v(g)\text{, \ where \ }N=|G|.
\mytag{2.1}
$$
We used the symbol $M$ for denoting the operation of invariant
averaging \mythetag{2.1} since it is analogous to the mathematical
expectation or the mean value in the theory of probability.\par
     Note that the operation of invariant averaging \mythetag{2.1}
can be applied not only to numeric functions, but to vector-valued,
operator-valued, and matrix-valued functions on a group $G$. In 
order to apply this operation to a function it should take its
values in some linear vector space. Then the result of averaging 
it $M[v]$ is an element of the same linear vector space. The operation
of invariant averaging satisfies the following obvious conditions
of linearity:
\roster
\item\quad $M[u+v]=M[u]+M[v]$;
\item\quad $M[\alpha\,u]=\alpha\,M[u]$, where $\alpha$ is a number.
\endroster
The invariance of the averaging \mythetag{2.1} reveals in the form of the
following relationships:
\roster
\item[3]\quad\kern -1.2cm\vtop{\hsize=9.6cm\noindent $M[R(g)u]=M[u]$, the 
invariance with respect to right shifts;}
\vskip 0.5ex
\item\quad\kern -1.2cm\vtop{\hsize=9.6cm\noindent $M[L(g)u]=M[u]$, the 
invariance with respect to left shifts;}
\vskip 0.5ex
\item\quad\kern -1.2cm\vtop{\hsize=9.6cm\noindent $M[Au]=M[u]$, the
invariance with respect to the inversion.}
\vskip 0.5ex
\endroster
The proof of the properties \therosteritem{3}--\therosteritem{5} is
reduced to verifying the following relationships:
$$
\align
M[R(g)u]&=\frac{1}{N}\sum_{a\in G}u(a\,g)=\frac{1}{N}\sum_{b\in G}
u(b)=M[u],\\
M[L(g)u]&=\frac{1}{N}\sum_{a\in G}u(g^{-1}\,a)=\frac{1}{N}\sum_{b\in G}
u(b)=M[u],\\
M[Au]&=\frac{1}{N}\sum_{a\in G}u(a^{-1})=\frac{1}{N}\sum_{b\in G}
u(b)=M[u].
\endalign
$$
Remember that the inversion operator $A$ in the property \therosteritem{5}
above is defined by the relationship \mythetag{1.4}.\par
     If $u$ is a vector-valued function with the values in a linear vector
space $V$ , then the properties \therosteritem{1}--\therosteritem{5} can be
complemented  with one more property:
\roster
\item[6]\quad\kern -1.2cm\vtop{\hsize=9.6cm\noindent $M[Bu]=BM[u]$, where 
$B$ is an arbitrary linear mapping with the domain $V$.}
\endroster
The relationship like \therosteritem{6} is fulfilled for operator-valued 
functions with the values in $\End(V)$:
\roster
\item[6]\quad\kern -1.2cm\vtop{\hsize=9.6cm\noindent $M[B\compos u]
=B\compos M[u]$, where $B$ is an arbitrary linear mapping with the 
domain $V$.}
\endroster
Moreover, for such functions with values in $\End(V)$ we can write the
following two additional properties:
\roster
\item[7]\quad $\tr M[u]=M[\tr u]$;
\item\quad\kern -1.2cm\vtop{\hsize=9.6cm\noindent $M[B\compos u]
=B\compos M[u]$, where $B$ is an arbitrary linear mapping with the 
domain $V$.}
\endroster
The operation of invariant averaging \mythetag{2.1} plays an important 
role in the theory of representations of finite groups. As the first 
example of its usage we prove the following fact.
\mytheorem{2.1} Each finite-dimensional representation of a finite group
is equivalent to some unitary representation of it.
\endproclaim
     Let $(f,G,V)$ be some finite-dimensional representation of a finite
group $G$. Generally speaking, in order to prove the theorem we should 
construct a unitary representation $(h,G,W)$  of the same group in some
Hermitian space $W$ and find a linear mapping $A\!:\,V\to W$ being an
isomorphism of representations $f$ and $h$. Assume for a while that we
managed to do it. Then we have the following relationships:
$$
\xalignat 2
&A\compos f(g)=h(g)\compos A,
&&\langle h(g)\bold u|h(g)\bold v\rangle=\langle\bold u|\bold v\rangle.
\endxalignat
$$
The space $V$ is not equipped with its own scalar product. However, we
equip it with a scalar product as follows:
$$
\hskip -2em
\langle\bold u|\bold v\rangle=\langle A\bold u|A\bold v\rangle.
\mytag{2.2}
$$
All properties of a scalar product for the sesquilinear form 
\mythetag{2.2} are verified immediately. The positivity is 
present because $A$ is a bijective mapping and $\Ker A=\{0\}$.
The representation $f$ appears to be a unitary representation 
with respect to the Hermitian scalar product \mythetag{2.2}:
$$
\gather
\langle f(g)\bold u|f(g)\bold v\rangle=\langle Af(g)\bold u|Af(g)\bold v
\rangle=\\
=\langle h(g)A\bold u|h(g)A\bold v\rangle=
\langle A\bold u|A\bold v\rangle=\langle\bold u|\bold v\rangle,
\endgather
$$
while $A$ establishes the unitary equivalence for $f$ and $h$. These
considerations show that in order to prove the theorem~\mythetheorem{2.1}
there is no need to construct a separate unitary representation 
$h,G,W$ and find an isomorphism $A$. It is sufficient to define a proper
scalar product in $V$ such that $f$ is a unitary representation with
respect to it.\par
   Let $\langle\!\langle f(g)\bold u|f(g)\bold v\rangle\!\rangle$ be some
arbitrary scalar product in $V$. For instance, it can be defined using
the coordinates of vectors $\bold u$ and $\bold v$ in some fixed basis:
$$
\langle\!\langle f(g)\bold u|f(g)\bold v\rangle\!\rangle
=\sum^n_{i=1}\overline{u^i}\,v^i.
$$
The operators $f(g)$ should not be unitary operators with respect to
such a scalar product. So, we need to improve it. Let's define another
scalar product in $V$ by means of the operation of invariant averaging:
$$
\langle\bold u|\bold v\rangle=
M[\langle\!\langle f(g)\bold u|f(g)\bold v\rangle\!\rangle]
=\frac{1}{N}\sum_{g\in G}\langle\!\langle f(g)\bold u|f(g)
\bold v\rangle\!\rangle.\quad
\mytag{2.3}
$$
It is easy to see that the form \mythetag{2.3} is sesquilinear and
symmetric. It is also a positive form:
$$
\langle\bold u|\bold u\rangle=
\sum_{g\in G}\frac{\langle\!\langle f(g)\bold u|f(g)
\bold u\rangle\!\rangle}{N}=\sum_{g\in G}\frac{\Vert f(g)\bold u\Vert^2}
{N}>0\text{\ \ for all \ }\bold u\neq 0.
$$
The operators $f(g)$ are unitary operators with respect to 
the scalar product \mythetag{2.3}. This fact follows from the 
property~\therosteritem{3} of the invariant averaging. Indeed, 
we have
$$
\gather
\langle f(g)\bold u| f(g)\bold v\rangle=\sum_{a\in G}\frac{\langle
\!\langle f(a)f(g)\bold u|f(a)f(g)\bold v\rangle\!\rangle}{N}=\\
=\sum_{a\in G}\frac{\langle\!\langle f(a\,g)\bold u|f(a\,g)\bold v
\rangle\!\rangle}{N}=\sum_{b\in G}\frac{\langle\!\langle f(b)\bold u
|f(b)\bold v\rangle\!\rangle}{N}=\langle\bold u|\bold v\rangle.
\endgather
$$
The above considerations prove that each finite-dimensional representation
of a finite group can be transformed to a unitary representation by means
of the proper choice \mythetag{2.3} of a scalar product. Thus, the
theorem~\mythetheorem{2.1} is proved.\par
     As an immediate corollary of the theorem~\mythetheorem{2.1} we 
get the following important proposition concerning finite dimensional
representations of finite groups.\par
\mytheorem{2.2} Each finite-dimensional representation of a finite
group is completely reducible.
\endproclaim
    The proof of this theorem is based on the
theorem~\mythetheoremchapter{7.5}{1} from
Chapter~\uppercase\expandafter{\romannumeral 1}. This theorem says
that each unitary representation is completely reducible. As for the
finite-dimensional representations of finite groups, we have proved
their equivalence to some unitary representations. 
\head
\SectionNum{3}{50} Characters of group representations.
\endhead
\rightheadtext{\S\,3. Characters of group representations.}
     Let $(f,G,V)$ be some finite-dimensional representation of a
group $G$. Each such representation is associated with the numeric
function $\chi_f$ on the group $G$ defined through the traces of 
representation operators:
$$
\hskip -2em
\chi_f(g)=\tr f(g).
\mytag{3.1}
$$
The numeric function $\chi_f(g)$ on $G$ introduced by the formula
\mythetag{3.1} is called the {\it character\/} of the representation
$f$.
\mytheorem{3.1} The characters of finite-dimensional representations
possess the following properties:
\roster
\rosteritemwd=2pt
\item the characters of two equivalent representations do coincide;
\item a character is constant within each conjugacy class;
\item if $f$ is a unitary representation, then $\chi_f(g^{-1})
      =\overline{\chi_f(g)}$;
\item the character of the direct sum of representations is equal to the
      sum of characters of separate direct summands;
\item the character of the tensor product of representations is equal to
      the product of characters of its multiplicands.
\endroster
\endproclaim
     Let's begin with proving the first item of the theorem. Assume that
we have two equivalent representation $(f,G,V)$ and $(h,G,W)$ and let 
$A\!:\,V\to W$ be an isomorphism of these representations. Let's choose
a basis $\bold e_1,\,\ldots,\,\bold e_n$ in $V$. Then the vectors
$\tilde\bold e_1=A\bold e_1,\,\ldots,\,\tilde\bold e_n=A\bold e_n$
constitute a basis in the space $W$. Let's calculate the matrices of the
operators $f(g)$ and $h(g)$ in these bases. They are defined by the
following relationships:
$$
\xalignat 2
&f(g)\bold e_i=\sum^n_{j=1}F^j_i(g)\,\bold e_j,
&&h(g)\tilde\bold e_i=\sum^n_{j=1}H^j_i(g)\,\tilde\bold e_j.
\qquad
\mytag{3.2}
\endxalignat
$$
Let's substitute $\tilde\bold e_i=A\bold e_i$ and $\tilde\bold e_j
=A\bold e_j$ into the second formula \mythetag{3.2} and take into
account the relationship $A\compos f(g)=h(g)\compos A$. As a result 
we obtain
$$
\hskip -2em
Af(g)\bold e_i=\sum^n_{j=1}H^j_i(g)\,A\bold e_j.
\mytag{3.3}
$$
The mapping $A$ is a bijective mapping. Therefore, we can cancel it in
\mythetag{3.3}. Upon canceling $A$, we compare \mythetag{3.3} with the
first formula \mythetag{3.2}. This comparison yields $H^j_i(g)=F^j_i(g)$,
i\.\,e\. the matrices of the operators $f(g)$ and $h(g)$ do coincide.
Hence, $\tr f(g)=\tr h(g)$ and $\chi_f(g)=\chi_h(g)$. The first item
in the theorem~\mythetheorem{3.1} is proved.\par
     Assume that $g$ and $\tilde g$ are two elements of the same conjugacy
class in $G$. Then $\tilde g=a\,g\,a^{-1}$ for some $a\in G$. Therefore,
we get
$$
f(\tilde g)=f(a\,g\,a^{-1})=f(a)\compos f(g)\compos f(a)^{-1}.
$$
Now it is sufficient to apply the formula $\tr(B\compos A \compos B^{-1})
=\tr(A)$ setting $A=f(g)$ and $B=f(a)$ in it. The second
item in the theorem~\mythetheorem{3.1} is proved.\par
     In order to prove the third item we consider a unitary representation
$(f,G,V)$ and choose some orthonormal basis in $V$. The condition that $f$
is unitary is written as $\langle f(g)\bold u|\bold v\rangle=\langle
\bold u|f(g)^{-1}\bold v\rangle$. Upon substituting $\bold u=\bold e_i$
and $\bold v=\bold e_j$ we take into account the first formula \mythetag{3.2}. 
As a result we get
$$
\hskip -2em
\overline{F^j_i(g)}=F^i_j(g^{-1}).
\mytag{3.4}
$$
The relationship \mythetag{3.4} means that the matrices $\overline{F(g)}$
and $F(g^{-1})$ are obtained from each other by transposing. The traces of 
any two such matrices do coincide. The third item is proved.\par
     The fourth item is trivial. Let $f=f_1\oplus f_2$ and $V=V_1\oplus
V_2$. We choose a basis in $V$ composed by two bases in $V_1$ and $V_2$
respectively. The matrix of the operator $f(g)$ in such a basis is a
blockwise diagonal matrix, the diagonal blocks of it coincides with the
matrices of the operators $f_1(g)$ and $f_2(g)$. Therefore,
$\tr f(g)=\tr f_1(g)+\tr f_2(g)$.\par
     Now let's proceed to the fifth item of the theorem. Let's denote
$\varphi=f\otimes h$ and $V=U\otimes W$. We choose some basis $\bold e_1,
\,\ldots,\,\bold e_n$ in $U$ and some basis $\tilde\bold e_1,
\,\ldots,\,\tilde\bold e_m$ in $W$. The matrices of the operators $f(g)$
and $h(g)$ are determined by the relationships
$$
\xalignat 2
&f(g)\bold e_i=\sum^n_{j=1}F^j_i(g)\,\bold e_j,
&&h(g)\tilde\bold e_q=\sum^n_{p=1}H^p_q(g)\,\tilde\bold e_p,
\qquad
\mytag{3.5}
\endxalignat
$$
which are analogous to \mythetag{3.2}. The vectors $\bold E_{iq}=
\bold e_i\otimes\tilde\bold e_q$ constitute a basis in the tensor
product $V=U\otimes W$. The vectors of this basis are enumerated 
by two indices, therefore the matrix of the operator $\varphi(g)$ 
in this basis is represented by a four-index array. It is determined 
by the following relationships:
$$
\hskip -2em
\varphi(g)\bold E_{iq}=\sum^n_{j=1}\sum^m_{p=1}\varPhi^{jp}_{iq}(g)
\,\bold E_{jp}.
\mytag{3.6}
$$
The action of the operator $\varphi(g)$ upon the basis vectors 
$\bold E_{iq}=\bold e_i\otimes\tilde\bold e_q$ is determined by the
formula \mythetagchapter{5.1}{1} from
Chapter~\uppercase\expandafter{\romannumeral 1}:
$$
\hskip -2em
\varphi(g)(\bold e_i\otimes\tilde\bold e_q)=(f(g)\bold e_i)\otimes
(h(g)\tilde\bold e_q).
\mytag{3.7}
$$
Combining \mythetag{3.7} with the formulas \mythetag{3.5}, we find
$$
\hskip -2em
\varphi(g)\bold  E_{iq}=\sum^n_{j=1}\sum^m_{p=1} F^j_i(g)\,
H^p_q(g)\,\bold  E_{jp}.
\mytag{3.8}
$$
Now, comparing the relationships \mythetag{3.6} and \mythetag{3.8}, we
determine the matrix components for the operator $\varphi(g)$: 
$$
\hskip -2em
\varPhi^{jp}_{iq}(g)=F^j_i(g)\,H^p_q(g).
\mytag{3.9}
$$
The rest is to calculate the trace of the operator $\varphi(g)$ as the
trace of a matrix array \mythetag{3.9}:
$$
\gather
\tr\varphi(g)=\sum^n_{i=1}\sum^m_{q=1}\varPhi^{iq}_{iq}=
\sum^n_{i=1}\sum^m_{q=1}F^i_i(g)\,H^q_q(g)=\\
\left(\,\shave{\sum^n_{i=1}}F^i_i(g)\right)
\left(\,\shave{\sum^m_{q=1}}H^q_q(g)\right)=
\tr f(g)\,\tr h(g).
\endgather
$$
This relationship completes the proof of the fifth item 
in the theorem~\mythetheorem{3.1} and the proof of the theorem 
in whole.\par
     In the end of this section one should remark that the properties
of representation character considered above are valid for 
finite-dimensional representations of arbitrary groups, not for 
finite groups only.
\head
\SectionNum{4}{54} Orthogonality relationships.
\endhead
\rightheadtext{\S\,4. Orthogonality relationships.}
     Let $(f,G,V)$ and $(h,G,W)$ are two complex finite-dimensional
representations of a finite group $G$. We choose some linear mapping 
$B\!:\,V\to W$ and, using it, we define a function $\varphi_B(g)$ on
$G$ with the values in $\Hom(V,W)$. Let's set
$$
\varphi_B(g)=h(g)\compos B\compos f(g^{-1}).
$$
The result of invariant averaging of $\varphi_B(g)$ over the group $G$
is some element $C\in\Hom(V,W)$:
$$
\hskip -2em
C=M[\varphi_B(g)]=\frac{1}{N}\sum_{a\in G}
h(a)\compos B\compos f(a^{-1}).
\mytag{4.1}
$$ 
It is easy to verify that the mapping $C\!:\,V\to W$ is a homomorphism
of the representations $f$ and $h$. Indeed, we have
$$
\gather
C\compos f(g)=M[\varphi_B(g)]\compos f(g)=\frac{1}{N}\sum_{a\in G}
h(a)\compos B\compos f(a^{-1})\compos f(g)=\\
=\frac{1}{N}\sum_{a\in G}h(a)\compos B\compos f(a^{-1}\,g)=
\frac{1}{N}\sum_{b\in G}h(g\,b)\compos B\compos f(b^{-1})=\\
=\frac{1}{N}\sum_{b\in G}h(g)\compos h(b)\compos B\compos f(b^{-1})=
h(g)\compos M[\varphi_B(g)]=h(g)\compos C.
\endgather
$$\par
      Assume that the representations $(f,G,V)$ and $(h,G,W)$ 
are irreducible. If they are not equivalent, applying Schur's 
lemma~\mythelemmachapter{5.1}{1}, we get $C=0$.\par
      Let's study the case $f\cong h$. For each pair of equivalent
irreducible representations we fix some bijective mapping $A_{fh}\!:
\,V\to W$ implementing an isomorphism of these representations. Then
the following lemma determines the structure of the mapping $C$ in
\mythetag{4.1}.
\mylemma{4.1} A homomorphism $C\!:\,V\to W$ of two equivalent irreducible
finite-dimensional complex representations $(f,G,V)$ and $(h,G,W)$ is fixed
up to a numeric factor, i\.\,e\. $C=\lambda\,A_{fh}$.
\endproclaim
\demo{Proof} Let's consider the operator $A=A_{fh}^{-1}\compos C$ in the
space $V$. Being the composition of two homomorphisms, this operator
implements an isomorphism of $f$ with itself. Therefore we have 
$A\compos f(g)=f(g)\compos A$ for all $f(g)$. Applying Schur's 
lemma~\mythelemmachapter{5.2}{1}, we get $A=\lambda\cdot 1$. Hence,
$C=\lambda\,A_{fh}$. The lemma is proved.
\qed\enddemo
     In order to calculate the numeric factor $\lambda$ we use the 
trace of the operator $A$, which is its numeric invariant:
$$
\lambda=\frac{\tr A}{\tr 1}=\frac{\tr A}{\dim V}=\frac{1}{\dim V}
\,\tr(A_{fh}^{-1}\compos C).
$$
Let's substitute the expression \mythetag{4.1} for $C$ into this formula:
$$
\gather
\lambda=\frac{1}{N\,\dim V}\sum_{a\in G}\tr(A_{fh}^{-1}\compos
h(a)\compos B\compos f(a^{-1}))=\\
=\frac{1}{N\,\dim V}\sum_{a\in G}\tr(f(a)\compos A_{fh}^{-1}\compos
B\compos f(a^{-1}))=\frac{\tr(A_{fh}^{-1}\compos B)}
{\dim V}.
\endgather
$$
Here again we used the formula $\tr(F\compos D\compos F^{-1})=\tr(D)$
with $F=f(a)$ and $D=A_{fh}^{-1}\compos B$. The result of calculating
the parameter $\lambda$ enables us to formula the following proposition.
\mytheorem{4.1} For arbitrary two irreducible finite-dimensional complex
representations $(f,G,V)$ and $(h,G,W)$ of a finite group $G$ the 
relationship
$$
\pagebreak
\sum_{a\in G}\frac{h(a)\compos B\compos f(a^{-1})}{N}
=\cases\qquad 0&\text{\ \ for \ }f\not\cong h;\\
\vspace{2ex}
\dfrac{\tr(A_{fh}^{-1}\compos B)}
{\dim V}&\text{\ \ for \ }f\not\cong h
\endcases\quad
\mytag{4.2}
$$
is fulfilled. It is valid for an arbitrary choice of a linear mapping
$B\!:\,V\to W$ in $\Hom(V,W)$.
\endproclaim
    The relationship \mythetag{4.2} is the basic orthogonality 
relationship in the theory of representations of finite groups. 
Let's consider the matrix form of this relationship. Assume that the
bases $\bold e_1,\,\ldots,\,\bold e_n$ and $\tilde\bold e_1,\,\ldots,
\,\tilde\bold e_m$ in the spaces $V$ and $W$ are chosen. They determine
the matrices $F^p_i(a)$ and $H^j_q(a)$ for the operators $f(a)$ and
$h(a)$ respectively. They also determine the matrix $B^q_p$ for the
mapping $B\in\Hom(V,W)$. In the case $f\not\cong h$ the bases in $V$ 
and $W$ are not related to each other. In the case $f\cong h$ it is
convenient to choose the first basis arbitrarily and then define the 
other by means of the relationship
$$
\hskip -2em
\tilde\bold e_i=A_{fh}\bold e_i,\quad i=1,\,\ldots,\,n.
\mytag{4.3}
$$
Under such a choice of bases the mapping $A_{fh}$ is represented by 
the unit matrix, while the matrices of the operators $f(a)$ and $h(a)$
do coincide: $F^p_i(a)=H^p_i(a)$. As for the mapping $B$, we choose it
so that the only nonzero element in its matrix $B^q_p=1$ is placed in 
the crossing of the $q$-th row and the $p$-th column. If all these
provisions are made, then the orthogonality relationship \mythetag{4.2} 
is rewritten as follows:
$$
\hskip -2em
\frac{1}{N}\sum_{a\in G}H^j_q(a)\,F^p_i(a^{-1})
=\cases\quad 0 &\text{for \ }f\not\cong h,\\
\vspace{1ex}
\dfrac{\delta^j_i\,\delta^p_q}{n}&\text{for \ }f\cong h.
\endcases
\mytag{4.4}
$$\par
    Assume that the representations $f$ and $h$ are unitary ones. 
Above we have already proved that any finite-dimensional complex
representation of a finite group can be replaced by some unitary
representation equivalent to it. And if $f$ and $h$ are equivalent, 
then they are unitary equivalent as well. For this reason we can 
assume that the mapping $A_{fh}$ is an isometry, while the bases
$\bold e_1,\,\ldots,\,\bold e_n$ and $\tilde\bold e_1,\,\ldots,
\,\tilde\bold e_m$ are orthonormal bases. Then the relationship 
\mythetag{4.4} is rewritten as follows: 
$$
\hskip -2em
\frac{1}{N}\sum_{a\in G}H^j_q(a)\,\overline{F^i_p(a)}
=\cases\quad 0 &\text{for \ }f\not\cong h,\\
\vspace{1ex}
\dfrac{\delta^j_i\,\delta^p_q}{n}&\text{for \ }f\cong h.
\endcases
\mytag{4.5}
$$
Note that the equality \mythetag{4.3} is compatible with the
orthonormality of the bases $\bold e_1,\,\ldots,\,\bold e_n$ 
and $\tilde\bold e_1,\,\ldots,\,\tilde\bold e_m$ since $A_{fh}$
is an isometry. In writing \mythetag{4.5} we used the relationship 
$$
F^p_i(a^{-1})=\overline{F^i_p(a)}
$$
because the matrices of unitary operators $f(a)$ in an orthonormal
basis are unitary matrices.\par
     Let's set $q=j$ and $p=i$ in the formula \mythetag{4.5} and
then sum over the indices $i$ and $j$. As a result we derive from
\mythetag{4.5} the following relationship for the characters of 
irreducible representations $f$ and $h$ of a finite group:
$$
\hskip -2em
\frac{1}{N}\sum_{a\in G}\tr(h(a))\,\overline{\tr(f(a))}
=\cases 0 &\text{for \ }f\not\cong h,\\
1 &\text{for \ }f\cong h.
\endcases
\mytag{4.6}
$$
Note that the representations $f$ and $h$ in \mythetag{4.6} are 
not unitary ones. The matter is that the characters of equivalent
representations do coincide, while $f$ and $h$, according to the 
theorem~\mythetheorem{2.1}, are equivalent to some unitary 
representations.
\mytheorem{4.2} The characters of two non-equivalent irredu\-cible
finite-dimensional complex representations of a finite group $G$
are orthogonal as the elements of the space $L_2(G)$.
\endproclaim
     The relationship \mythetag{4.6} is a proof of the 
theorem~\mythetheorem{4.2}. In order to see it one should compare
this relationship with \mythetag{1.1}. From the finiteness
$\dim L_2(G)=N\leqslant\infty$ we conclude that the number of 
non-equivalent irreducible finite-dimensional complex 
representations of a finite group $G$ is finite. Therefore, 
considering a {\it complete set\/} of such representations 
is a true idea.
\mydefinition{4.1} The representations $f_1,\,\ldots,\,f_m$ form
a complete set of non-equivalent irreducible finite-dimensional
complex representations of a finite group $G$ if
\roster
\rosteritemwd=5pt
\item any two of them are not equivalent to each other;
\item each irreducible finite-dimensional complex representation
of $G$ is equivalent to one of the representations $f_1,\,\ldots,
\,f_m$.
\endroster
\enddefinition
     The number $m$ of the representations in a complete set is a 
numeric invariant of a finite group $G$. It is not greater than 
the order of the group $N=|G|$.\par
    Let $(f_1,G,V_1),\,\ldots,\,(f_m,G,V_m)$ be a complete set of 
non-equivalent irreducible representations. Without loss of generality
we can assume these representations to be unitary ones. Let $n_1,\,
\ldots,\,n_m$ be the dimensions of these representations. We choose
some orthonormal basis in each of the spaces $V_1,\,\ldots,\,V_m$.
Then we have a set of matrices with the components
$$
F^j_i(g,r),\quad r=1,\,\ldots,\,m;\quad 1\leqslant i,j\leqslant n_r.
$$
Each component in these matrices depend on $g\in G$, therefore, it can
be treated as a function from $L_2(G)$. From \mythetag{4.5} we derive 
the following orthogonality relationships for these functions:
$$
\hskip -2em
\frac{1}{N}\sum_{a\in G}F^j_q(a,r)\,\overline{F^i_p(a,s)}
=\frac{1}{n_r}\,\delta_{rs}\,\delta_{ij}\,\delta_{p\kern 0.5pt q}.
\mytag{4.7}
$$
The relationships \mythetag{4.7} mean that the functions $F^j_i(g,r)$
treated as the elements of the space $L_2(G)$ are pairwise orthogonal
to each other. Apparently, they are not only orthogonal, but form a
complete orthogonal set of functions in this space.
\mytheorem{4.3} For an arbitrary complete set $f_1,\,\ldots,\,f_m$ of
irreducible unitary representations of a finite group $G$ the matrix
elements of the operators $f_r(g)$ calculated in some orthonormal 
bases form a complete orthogonal system of functions in $L_2(G)$.
\endproclaim
\demo{Proof} The orthogonality of the functions $F^i_j(g,r)$ follows 
from the relationship \mythetag{4.7}. We need to prove their completeness.
For this purpose we consider the right regular representation
$(R,G,L_2(G))$. It is unitary with respect to the Hermitian structure
given by the scalar product \mythetag{1.1} (see 
theorem~\mythetheorem{1.1}). For this reason the right regular
representation $(R,G,L_2(G))$ is completely reducible, it is expanded 
into the direct sum of unitary irreducible representations:
$$
\hskip -2em
R=R_1\oplus\ldots\oplus R_k.
\mytag{4.8}
$$
The expansion \mythetag{4.8} of the right regular representation $R$
is associated with the expansion of the space $L_2(G)$ into the direct
sum of irreducible $R$-invariant subspaces 
$$
L_2(G)=W_1\oplus\ldots\oplus W_k.
$$
Each of the irreducible representations $R_q$ in \mythetag{4.8} is
equivalent to one of the irreducible unitary representations 
$(f_{r(q)},G,V_{r(q)})$ from our complete set $f_1,\,\ldots,\,f_m$.
Applying the theorem~\mythetheoremchapter{7.6}{1} from
Chapter~~\uppercase\expandafter{\romannumeral 1}, we conclude that
the representations $R_q$ and $f_{r(q)}$ are unitary equivalent.
Therefore, in each subspace $W_q\subseteq L_2(G)$ we can choose 
some orthonormal basis of functions
$$
\hskip -2em
\varphi_i(g,q),\quad 1\leqslant i\leqslant n_{r(q)},
\mytag{4.9}
$$
such that the matrices of the operators $R_q(g)$ coincide with the
matrix components $F^j_i(g,r(q))$ of the operators for the corresponding
representation $f_{r(q)}$ in the complete set. Let's write this fact
as a formula:
$$
\hskip -2em
R_q(g)\varphi_i(q,q)=\sum^{n_{r(q)}}_{j=1}F^j_i(q,r(q))\,\varphi_j(a,q).
\mytag{4.10}
$$
But $R_q(g)$ is the restriction of the operator $R(q)$ from \mythetag{4.8}
to its invariant subspace $W_q$, \pagebreak while $\varphi_i(a,q)$ is an 
element of this subspace. Therefore, we have
$$
\hskip -2em
R_q(g)\varphi_i(a,q)=R(g)\varphi_i(a,q)=\varphi_i(a\,g,q).
\mytag{4.11}
$$
Let's substitute \mythetag{4.11} into \mythetag{4.10}. Then in the
relationship obtained we set $a=e$. The quantities $\varphi_i(e,q)$
are some constant numbers, we denote them $c_{iq}=\varphi_i(e,q)$.
As a result we get
$$
\hskip -2em
\varphi_i(g,q)=\sum^{n_{r(q)}}_{j=1}c_{jq}\,F^j_i(q,r(q)).
\mytag{4.12}
$$
The formula \mythetag{4.12} is an expansion of the function 
$\varphi_i(g,q)$ in the set of functions $F^j_i(q,r(q))$. But the
set of functions $\varphi_i(g,q)$ in \mythetag{4.9} constitute a basis 
in $L_2(G)$. It is a complete set, each element $\varphi(g)\in L_2(G)$ 
has an expansion in this set of functions. Due to the formula \mythetag{4.12} 
such an expansion can be transformed into the expansion of $\varphi(g)$
in the set of functions $F^j_i(q,r(q))$. Hence, the set of functions
$F^j_i(q,r(q))$ also is a complete set for $L_2(G)$. The 
theorem~\mythetheorem{4.3} is proved.
\qed\enddemo
     Let $(\varphi,G,V)$ be some finite-dimensional complex 
representation of a finite group $G$. It is completely reducible 
(see theorem~\mythetheorem{2.2}). It is expanded into
a sum of irreducible representations:
$$
\hskip -2em
\varphi=\varphi_1\oplus\ldots\oplus\varphi_\nu.
\mytag{4.13}
$$
Each of the irreducible representation $\varphi_q$ in \mythetag{4.13}
is equivalent to one of the representations $f_{r(q)}$ in our complete 
set $f_1,\,\ldots,\,f_m$. Let's denote by $k_r$ the number of irreducible
representations in the expansion \mythetag{4.13} which are equivalent to
the representation $f_r$. Then the expansion \mythetag{4.13} is rewritten
as
$$
\hskip -2em
\varphi\cong k_1\,f_1\oplus\ldots\oplus k_m\,f_m.
\mytag{4.14}
$$
The number $k_r$ in \mythetag{4.14} is called the {\it multiplicity}
of the entry of the irreducible representation $f_r$ in $\varphi$. 
The orthogonality relationship \mythetag{4.6} for characters enables 
us to calculate the multiplicities without performing the expansion
\mythetag{4.13} itself:
$$
\hskip -2em
k_r=\frac{1}{N}\sum_{g\in G}\tr\varphi(g)\,\overline{\tr f_r(g)}.
\mytag{4.15}
$$
The relationship \mythetag{4.15} is derived from the following expansion
for the function $\tr\varphi(g)$:
$$
\hskip -2em
\tr\varphi(g)=k_1\,\tr f_1+\ldots+k_m\,\tr f_m.
\mytag{4.16}
$$
The relationship \mythetag{4.16} in turn is derived from \mythetag{4.14}.
\par
     Let's find the expansion \mythetag{4.14} in the case of the right
regular representation $R(g)$. The corresponding expansion for the left
regular representation is the same because these two representations are
equivalent. In order to calculate the trace of the operator $R(g)$ it
is necessary to choose a basis in $L_2(G)$ and calculate the matrix
elements of the operator $R(g)$ in this basis. The
theorem~\mythetheorem{4.3}, which was proved above, says that the set of
functions $F^j_i(g,r)$ corresponding to some complete set of irreducible
representations can be used as a basis in $L_2(G)$. The basis functions
$F^j_i(g,r)$ are enumerated with three indices $i$, $j$, and $r$. 
Therefore, the matrix of the operator $R(g)$ in this basis is represented
by a six-index array $R^{jp}_{qi}(r,s)$. This array is determined as
follows:
$$
\hskip -2em
R(g)F^j_i(a,r)=\sum^m_{s=1}\sum^{n_s}_{q=1}\sum^{n_s}_{p=1}
R^{jp}_{qi}(r,s)\,F^q_p(a,s).
\mytag{4.17}
$$
The trace of the six-index array $R^{jp}_{qi}(r,s)$ representing a matrix
is calculated according to the following formula:
$$
\hskip -2em
\tr R(g)=\sum^m_{r=1}\sum^{n_s}_{i=1}\sum^{n_s}_{j=1}
R^{ji}_{ji}(r,r).
\mytag{4.18}
$$
Let's calculate the left hand side of the equality \mythetag{4.17} 
directly from the relationship \mythetag{1.2} that define the operator
$R(g)$:
$$
R(g)F^j_i(a,r)=F^j_i(a\,g,r)=\sum^{n_r}_{p=1}F^j_p(a,r)\,F^p_i(g,r).
\quad
\mytag{4.19}
$$
Here, deriving the formula \mythetag{4.19}, we used the relationship
$f_r(a\,g)=f_r(a)\compos f_r(g)$ written in the matrix form. Comparing
the formulas \mythetag{4.17} and \mythetag{4.19}, we find
$$
R^{jp}_{qi}(r,s)=\delta_{rs}\delta^j_q\,F^p_i(g,r).
$$
The rest is to substitute this expression into \mythetag{4.18} and
perform the summations prescribed by that formula:
$$
\hskip -2em
\tr R(g)=\sum^m_{r=1}n_r\,\tr f_r(g).
\mytag{4.20}
$$
Let's compare \mythetag{4.20} and \mythetag{4.16}. The result of this
comparison is formulated in the following theorem.
\mytheorem{4.4} Each irreducible representation $f_r$ from some complete
set $f_1,\,\ldots,\,f_m$ of irreducible finite-dimensional complex
representations of a finite group $G$ enters the right regular
representation $(R,G,L_2(G))$ with the multiplicity $k_r$ equal to its
dimension, i\.\,e\. $k_r=n_r=\dim V_r$.
\endproclaim
     Exactly the same proposition is valid for the left regular 
representation $(L,G,L_2(G))$ of the group $G$ either. The 
theorem~\mythetheorem{4.4} has the following immediate corollary
that follows from the fact that $\dim L_2(G)=|G|$.
\mycorollary{4.1} The order of a finite group $N=|G|$ is equal to the
sum of squares of the dimensions of all its non-equivalent irreducible
finite-dimensional complex representations $f_1,\,\ldots,\,f_m$, 
i\.\,e\. $N=(n_1)^2+\ldots+(n_m)^2$.
\endproclaim
     The same result can be obtained if we calculate the total number 
of functions entering the orthogonality relationships \mythetag{4.7} and
forming a complete orthogonal set of functions in $L_2(G)$.\par
     Let's consider the set of characters $\chi_1,\,\ldots,\,\chi_m$ for
irreducible representations from some complete set $f_1,\,\ldots,\,f_m$.
Due to the theorem~\mythetheorem{4.2} and the relationships \mythetag{4.6} 
they are orthogonal. But in general case they do not form a complete set
of such functions in $L_2(G)$. From the theorem~\mythetheorem{3.1} we 
know that these functions are constants within conjugacy classes of $G$.
Let's denote by $M_2(G)$ the set of complex numeric functions on $G$
constant within each conjugacy class of $G$. It is clear that $M_2(G)$
is a linear subspace in $L_2(G)$. It inherits the Hermitian scalar product 
\mythetag{1.1} from the space $L_2(G)$.
\mytheorem{4.5} The characters $\chi_1,\,\ldots,\,\chi_m$ of the
representations $f_1,\,\ldots,\,f_m$ forming a complete set of irreducible
finite-dimensional representations of a finite group $G$ form a complete
set of orthogonal functions (a basis) in the space $M_2(G)\subseteq L_2(G)$.
\endproclaim
    Due to the theorem~\mythetheorem{3.1} all of the characters 
$\chi_1,\,\ldots,\,\chi_m$ belong to $M_2(G)$. They are orthogonal 
to each other and normalized to the unity. It follows from the 
theorem~\mythetheorem{4.2}. Let $\varphi(g)$ be some arbitrary element
of the space $M_2(G)$. Then we can expand it in the set of functions
$F^j_i(g,r)$ (see theorem~\mythetheorem{4.3}):
$$
\hskip -2em
\varphi(g)=\sum^m_{r=1}\sum^{n_r}_{i=1}\sum^{n_r}_{j=1}c_j^{\,i}(r)
\,F^j_i(g,r).
\mytag{4.21}
$$
Let's perform the conjugation $g\mapsto a\,g\,a^{-1}$ in the argument of
the function $\varphi(g)$. This operation does not change its value since
$\varphi(g)\in M_2(G)$. This value is not changed upon averaging over 
conjugations by means of all elements of the group $G$ either:
$$
\hskip -2em
\varphi(g)=\frac{1}{N}\sum_{a\in G}\varphi(a\,g\,a^{-1}).
\mytag{4.22}
$$
Let's substitute the expansion \mythetag{4.21} into \mythetag{4.22}. This
yields
$$
\varphi(g)=\sum^m_{r=1}\sum^{n_r}_{i=1}\sum^{n_r}_{j=1}c_j^{\,i}(r)
\left(\frac{1}{N}\shave{\sum_{a\in G}}F^j_i(a\,g\,a^{-1},r)\right).
\quad
\mytag{4.23}
$$
We denote by $\psi^j_i(g,r)$ the expression enclosed into round brackets
in right hand side of \mythetag{4.23}. For this quantity we get
$$
\psi^j_i(g,r)=\frac{1}{N}\sum_{a\in G}\sum^{n_r}_{p=1}\sum^{n_r}_{q=1}
F^j_p(a,r)\,F^p_q(g,r)\,F^q_i(a^{-1},r).
$$
Here we used the relationship $f_r(a\,g\,a^{-1})=f_r(a)\compos f_r(g)
\compos f_r(a^{-1})$ written in the matrix form. Now let's recall that
$F^j_i(a,r)$ are unitary matrices and apply the orthogonality relationship 
\mythetag{4.7}:
$$
\gather
\psi^j_i(g,r)=\frac{1}{N}\sum_{a\in G}\sum^{n_r}_{p=1}\sum^{n_r}_{q=1}
F^j_p(a,r)\,F^p_q(g,r)\,\overline{F^i_q(a,r)}=\\
=\frac{1}{n_r}\sum^{n_r}_{p=1}\sum^{n_r}_{q=1}F^p_q(g,r)
\,\delta^j_i\,\delta^q_p.
\endgather
$$
Upon performing the summation in the right hand side of the above formula
for the quantity $\psi^j_i(g,r)$ we get
$$
\psi^j_i(g,r)=\frac{1}{n_r}\tr f_r(g)\,\delta^j_i=\frac{1}{n_r}
\chi_r(g)\,\delta^j_i.
$$
The rest is to substitute this expression back into the formula
\mythetag{4.23}. As a result we find
$$
\hskip -2em
\varphi(g)=\sum^m_{r=1}\left(\,\shave{\sum^{n_r}_{i=1}}
\frac{c^{\,i}_i(r)}{n_r}\right)\chi_r(g).
\mytag{4.24}
$$
From \mythetag{4.24} it is clear that any function $\varphi(g)\in 
M_2(G)$ is expanded in the set of functions \pagebreak $\chi_1,\,
\ldots,\,\chi_m$.\par
    The following theorem on the number of irreducible representations 
in a complete set of such representations for a finite group is an obvious
corollary of the theorem~\mythetheorem{4.5}, which is proved just above. 
\mytheorem{4.6} The number of representations in a complete set 
$f_1,\,\ldots,\,f_m$ of irreducible finite-dimensional complex
representations of a finite group $G$ coincides with the number of
conjugacy classes in the group $G$.
\endproclaim
\head
\SectionNum{5}{65} Expansion into irreducible components.
\endhead
\rightheadtext{\S\,5. Expansion into irreducible components.}
     Let $G$ be a finite group. The theorem~\mythetheorem{4.6} determines
the number of irreducible representations in a complete set $f_1,\,\ldots,
\,f_m$, while the theorem~\mythetheorem{4.4} yields a way for finding such
representations. Indeed, each of the representations $f_1,\,\ldots,\,f_m$
enters the right regular representation $(R,G,L_2(G))$ at least once. 
Therefore, in order to find it one should construct the expansion 
\mythetag{4.13} for $\varphi=R$. For each particular finite group $G$
this could be done with the tools of linear algebra.\par
     Suppose that this part of work is already done and some complete set 
of unitary representations $f_1,\,\ldots,\,f_m$ is constructed. Then
upon choosing orthonormal bases in the spaces $V_1,\,\ldots,\,V_m$ of
these representations we can assume that the matrix elements $F^j_i(g,r)$
of the operators $f_r(g)$ are known. Under these assumptions we consider
the problem of expanding of a given representation $\varphi,G,V)$ into
its irreducible components. Let's begin with defining the following
operators:
$$
\hskip -2em
P^i_j(r)=\frac{n_r}{N}\sum_{a\in G}\overline{F^j_i(a,r)}\,\varphi(a).
\mytag{5.1}
$$
The number of such operators coincides with the number of functions  
$F^j_i(a,r)$. However, a part of these operators can be equal to zero. 
The operators $P^i_j(r)$ are interpreted as the coefficients for the
Fourier expansion of the operator-valued function in orthogonal system 
of functions in $L_2(G)$. The following expansion approves such 
interpretation:
$$
\hskip -2em
\varphi(g)=\sum^m_{r=1}\sum^{n_r}_{i=1}\sum^{n_r}_{j=1}
F^j_i(a,r)\,P^i_j(r).
\mytag{5.2}
$$
It is easily derived from the orthogonality relationship \mythetag{4.7}.
On the base of the same relationship \mythetag{4.7} one can derive a
number of other relationship for the operators \mythetag{5.1}. First of
all we consider the following ones:
$$
\align
&\hskip -2em
\varphi(g)\compos P^i_j(r)=\sum^{n_r}_{q=1}F^q_j(g,r)\,P^i_q(r),
\mytag{5.3}\\
&\hskip -2em
P^i_j(r)\compos\varphi(g)=\sum^{n_r}_{q=1}F^i_q(g,r)\,P^q_j(r).
\mytag{5.4}\\
\endalign
$$
We prove the relationship \mythetag{5.3} by means of direct calculations.
In order to transform the left hand side of \mythetag{5.3} we use the
formula \mythetag{5.1} for the operator $P^i_j(r)$:
$$
\gather
\varphi(g)\compos P^i_j(r)=\frac{n_r}{N}\sum_{a\in G}
\overline{F^j_i(a,r)}\,\varphi(g)\compos\varphi(a)=\\
=\frac{n_r}{N}\sum_{a\in G}\overline{F^j_i(a,r)}
\,\varphi(g\,a).
\endgather
$$
Replacing $a$ with $b=g\,a$ in summation over the group, we get
$$
\gather
\varphi(g)\compos P^i_j(r)=\frac{n_r}{N}\sum_{b\in G}
\overline{F^j_i(g^{-1}\,b,r)}\,\varphi(b)=\\
=\frac{n_r}{N}\sum_{b\in G}\sum^{n_r}_{q=1}\overline{
F^j_q(g^{-1},r)\,F^q_i(b,r)}\,\varphi(b).
\endgather
$$
Here we used the relationship $f_r(g^{-1}\,b)=f_r(g^{-1})\compos
f_r(b)$ written in the matrix form. Now the rest is to use the
relationship $f_r(g^{-1})=f_r(g)^{-1}$ and the unitarity of the
matrix $F(g,r)$:
$$
\gather
\varphi(g)\compos P^i_j(r)=\sum^{n_r}_{q=1}F^q_j(g,r)\left(
\frac{n_r}{N}\shave{\sum_{b\in G}}\overline{F^q_i(b,r)}\,\varphi(b)
\right)=\\
=\sum^{n_r}_{q=1}F^q_j(g,r)\,P^i_q(r).
\endgather
$$
Comparing the left and right hand sides of the above equalities, we 
see that the formula \mythetag{5.3} is proved. The formula \mythetag{5.4}
is proved in a similar way, therefore, we do not give its proof here.\par
     Let's set $j=i$ in the formulas \mythetag{5.3} and \mythetag{5.4}
and then sum up these equalities over the index $i$. The double sums in
right hand sides of the resulting equalities do coincide. Therefore, the
result can be written as
$$
\hskip -2em
\gathered
\sum^{n_r}_{i=1}\varphi(g)\compos P^i_i(r)=\sum^{n_r}_{i=1}
P^i_i(r)\compos\varphi(g)=\\
=\sum^{n_r}_{i=1}\sum^{n_r}_{j=1}F^j_i(g,r)\,P^i_j(r).
\endgathered
\mytag{5.5}
$$
The right hand side of \mythetag{5.5} differs from that of 
\mythetag{5.2} by the absence of the sum over $r$. Combining 
\mythetag{5.5} and \mythetag{5.2}, we obtain
$$
\hskip -2em
\varphi(g)=\sum^m_{r=1}\sum^{n_r}_{i=1}\varphi(g)\compos P^i_i(r)
\sum^m_{r=1}\sum^{n_r}_{i=1}P^i_i(r)\compos\varphi(g).
\mytag{5.6}
$$
Due to \mythetag{5.6} it is natural to introduce new operators by
means of the following formula:
$$
\hskip -2em
P(r)=\sum^{n_r}_{i=1}P^i_i(r)=\frac{n_r}{N}\sum_{a\in G}
\overline{\tr f_r(a)}\,\varphi(a).
\mytag{5.7}
$$
In terms of \mythetag{5.7} the relationship \mythetag{5.6} itself 
is rewritten as 
$$
\hskip -2em
\varphi(g)=\sum^m_{r=1}\varphi(g)\compos P(r)
\sum^m_{r=1}P(r)\compos\varphi(g).
\mytag{5.8}
$$
Setting $g=e$ in \mythetag{5.8}, we get an expansion of the unity (of the
identical operator) in operators \mythetag{5.7}:
$$
\hskip -2em
1=\sum^m_{r=1}P(r).
\mytag{5.9}
$$
\mytheorem{5.1} The operators $P(r)\!:\,V\to V$, $r=1,\,\ldots,\,m$, given
by the formula \mythetag{5.7} possess the following properties:
\roster
\rosteritemwd=5pt
\item they satisfy the relationships $P(r)^2=P(r)$, because of which 
those of them being nonzero $P(r)\neq 0$ are projectors onto the subspaces
$V(r)=\Img P(r)$;
\item they commute with the representation operators $\varphi(g)$, because
of which the subspaces $V(r)$ are invariant with respect to $\varphi(g)$;
\item they satisfy the relationship \mythetag{5.9} and the relationships
\linebreak
$P(r)\compos P(s)=0$ for $r\neq s$, because of which the expansion
$V=V(1)\oplus\ldots\oplus V(m)$ is an expansion into the direct sum of
invariant subspaces.
\endroster
\endproclaim
    The relationships $P(r)^2=P(r)$ from the first item of the theorem
and the relationships $P(r)\compos P(s)=0$ for $r\neq s$ from the third
item of the theorem can be combined into one relationship:
$$
\hskip -2em
P(r)\compos P(s)=P(r)\,\delta_{rs}=\cases 0 &\text{for \ } r\neq s,\\
P(r) &\text{for \ }r=s.\endcases
\mytag{5.10}
$$
The relationship \mythetag{5.10} is easily derived from the following more
general relationship for the operators $P^i_j(r)$ defined in \mythetag{5.1}:
$$
\hskip -2em
P^i_j(r)\compos P^k_q(s)=\delta_{rs}\,\delta^i_q\, P^k_j(r).
\mytag{5.11}
$$
It is convenient to prove \mythetag{5.11} by direct calculations. From the
formula \mythetag{5.1} for the product of operators in the left hand side
of \mythetag{5.11} we derive
$$
P^i_j(r)\compos P^k_q(s)=\frac{n_r\,n_s}{N^2}\sum_{a\in G}\sum_{b\in G}
\overline{F^j_i(a,r)\,F^q_k(b,s)}\,\varphi(a\,b).
$$
Let's denote $c=a\,b$ and choose $c$ as a new parameter in summation over
the group $G$ in place of $b$. This yields
$$
\gather
P^i_j(r)\compos P^k_q(s)=\frac{n_r\,n_s}{N^2}\sum_{a\in G}\sum_{c\in G}
\overline{F^j_i(a,r)\,F^q_k(a^{-1}\,c,s)}\,\varphi(c)=\\
=\frac{n_r\,n_s}{N^2}\sum_{a\in G}\sum_{c\in G}\overline{F^j_i(a,r)}
\sum^{n_s}_{p=1}\overline{F^q_p(a^{-1},s)\,F^p_k(c,s)}\,\varphi(c).
\endgather
$$
In the next step we use the unitarity of the matrix $F(a,s)$:
$$
P^i_j(r)\compos P^k_q(s)=\sum^{n_s}_{p=1}\sum_{a\in G}
\frac{n_r\,\overline{F^j_i(a,r)}\,F^p_q(a,s)}{N}
\sum_{c\in G}\frac{n_s\,\overline{F^p_k(c,s)}\,\varphi(c)}{N}.
$$
And finally, we use \mythetag{5.1} and the orthogonality relationship
\mythetag{4.7}:
$$
P^i_j(r)\compos P^k_q(s)=\sum^{n_s}_{p=1}\delta_{rs}\,\delta^p_j\,
\delta^i_q\,P^k_p(s)=\delta_{rs}\,\delta^i_q\,P^k_j(r).
$$
Comparing the left and right hand sides of this formula with
\mythetag{5.11}, we see that the formula \mythetag{5.11} is proved.
Hence, the relationship \mythetag{5.10} is also proved.\par
     The relationship $P(r)^2=P(r)$, which follows from \mythetag{5.10},
in the case of a nonzero operator $P(r)\neq 0$ means that $P(r)$ is a 
projection operator. It projects the space $V$ onto the subspace
$V(r)=\Img P(r)$ parallel to the subspace $\Ker P(r)$ (see more details
in \mybookcite{1}). The relationship 
$$
\hskip -2em
P(r)\compos P(s)=0=P(s)\compos P(r)\text{\ \ for \ }r\neq s
\mytag{5.12}
$$
means that the projectors $P(1),\,\ldots,\,P(m)$ commute and that 
$\Img P(r)\subseteq\Ker P(s)$ for $r\neq s$. Due to \mythetag{5.12} 
and \mythetag{5.9} the set of projectors $P(1),\,\ldots,\,P(m)$ is
a concordant and complete set of projectors. This set of projectors
determines an expansion of $V$ into a direct sum of subspaces:
$$
\hskip -2em
V=V(1)\oplus\ldots\oplus V(m)=\bigoplus^m_{r=1}V(r).
\mytag{5.13}
$$
The second item of the theorem~\mythetheorem{5.1} claiming the
commutativity of $P(r)$ and $\varphi(g)$ follows immediately from
\mythetag{5.5}. Due to this fact all subspaces in the expansion
\mythetag{5.13} are invariant subspaces of the representation
$(\varphi,G,V)$.\par
     The matrix elements $F^j_i(a,r)$ of the operators $f_r(a)$ 
depend on a choice of bases in spaces where they act. The operators
$P^i_j(r)$ calculated according to the formula \mythetag{5.1}
also depend on a choice of these bases. However, the operators
$P(r)$ do not depend on bases since the traces of the operators
$f_r(a)$ in the formula \mythetag{5.7} are their (basis-free)
scalar invariants. Therefore, the expansion \mythetag{5.13} is
also invariant, it is determined by the group $G$ itself and by its 
representation $\varphi$. Let's study how the expansions \mythetag{5.13}
and \mythetag{4.14} are related to each other.
\mytheorem{5.2} The restriction of the representation $\varphi$ to
the invariant subspace $V(r)=\Img P(r)$ is isomorphic to the irreducible
representation $f_r$ taken with the multiplicity $k_r$, where $k_r$ is
the coefficient of $f_r$ in the expansion \mythetag{4.14}.
\endproclaim
     In order to prove the theorem~\mythetheorem{5.2} we use the following
relationship whose left hand side coincides with the operator $P(r)$ in the
case of $\varphi=f_s$ (see formula \mythetag{5.7}):
$$
\hskip -2em
\frac{n_r}{N}\sum_{a\in G}\sum^{n_r}_{i=1}\overline{F^i_i(a,r)}
\,f_s(a)=\cases 0\text{\ \ for \ } s\neq r,\\ 1 \text{\ \ for \ }
s=r.\endcases
\mytag{5.14}
$$
In order to verify that \mythetag{5.14} \pagebreak is valid it 
is sufficient to pass from the operators $f_s(a)$ to their matrices 
$F(a,s)$ and then to apply the orthogonality relationship 
\mythetag{4.7}.\par
    Now let's consider the expansion \mythetag{4.13}. According to 
\mythetag{4.13}, the space $V$ is a direct sum of irreducible subspaces
$V=V_1\oplus\ldots\oplus V_\nu$, the restriction of $\varphi$ to each
such subspace is isomorphic to some irreducible representation from the
complete set $f_1,\,\ldots,\,f_r$:
$$
\hskip -2em
\varphi\,\hbox{\vrule height 8pt depth 10pt width 0.5pt}_{\,V_q}
\cong f_{r(q)}.
\mytag{5.15}
$$
Substituting \mythetag{5.15} into \mythetag{5.7} and taking into account 
\mythetag{5.14}, we find that $V(r)$ is the sum of those subspaces $V_q$
in the expansion $V=V_1\oplus\ldots\oplus V_\nu$ for which $r(q)=r$. The
number of such subspaces is equal to $k_r$, while the restriction of
$\varphi$ to each of them is isomorphic to $f_r$. The 
theorem~\mythetheorem{5.2} is proved.
     According to the theorem~\mythetheorem{5.2}, the operators $P(r)$
yield a constructive way to build the expansion \mythetag{4.14}, while
the expansion itself is unique up to a permutation of the summands.\par
     Let's consider a separate subspace $V(r)$ corresponding to the
component $k_r\,f_r$ in the expansion \mythetag{4.14}. If $k_r=0$, the
subspace $V(r)=\{0\}$ is trivial. If $k_r=1$ the subspace $V(r)$ is 
irreducible, it does not require a further expansion. The rest is the
case $k_r>1$, it should be especially treated. In this case the subspace
$V(r)$ is expanded into the sum of several irreducible subspaces
$$
\hskip -2em
V(r)=\bigoplus^{k_r}_{q=1}W_q\text{, \ where \ }\dim W_q=n_r.
\mytag{5.16}
$$
In contrast to the expansion \mythetag{5.13}, the expansion \mythetag{5.16}
is not unique. One of the ways for constructing such an expansion is due to
the operators $P^i_i(r)$. Their sum is equal to $P(r)$ according to the
formula \mythetag{5.7}. From \mythetag{5.11} for these operators we derive
$$
\xalignat 2
&P^i_i(r)^2=P^i_i(r),
&&P^i_i(r)\compos P^j_j(r)=0\text{\ \ for \ } i\neq j.\qquad\quad
\mytag{5.17}
\endxalignat
$$
There is no summation over $i$ and $j$ in these formulas. Due to 
\mythetag{5.17} the operators $P^i_i(r)$, $i=1,\,\ldots,\,n_r$ form 
a concordant and complete set of projection operators. They define 
an expansion of $V(r)$ into a direct sum of smaller subspaces:
$$
\hskip -2em
V(r)=\bigoplus^{n_r}_{i=1}V_i(r).
\mytag{5.18}
$$\par
     The projection operators $P^i_i(r)$ do not commute with 
$\varphi(g)$. Therefore, the subspaces $V_i(i)=\Img P^i_i(r)$ 
in the expansion \mythetag{5.18} are not invariant for the 
representation $\varphi$. However, we can overcome this difficulty.
For this purpose we use the following equalities easily derived from
\mythetag{5.11}:
$$
\xalignat 2
&P^k_k(r)\compos P^i_k(r)=P^i_k(r),
&&P^i_k(r)\compos P^i_i(r)=P^i_k(r),\qquad\quad
\mytag{5.19}\\
&P^k_i(r)\compos P^i_k(r)=P^i_i(r),
&&P^i_k(r)\compos P^k_i(r)=P^k_k(r).\qquad\quad
\mytag{5.20}
\endxalignat
$$
\mytheorem{5.3} For $i\neq k$ the operator $P^i_k(r)$ performs 
a bijective mapping from $V_i(r)$ onto $V_k(r)$.
\endproclaim
    The proof of the theorem~\mythetheorem{5.3} is based on the
formulas \mythetag{5.19} and \mythetag{5.20}. Indeed,  let $\bold v
\in V_i(r)$ and let $\bold u=P^i_k(r)\bold v$. Then, using the first
relationship \mythetag{5.19}, we derive $P^k_k\bold u=\bold u$. Hence, 
$\bold u\in V_k(r)$, i\.\,e\. $P^i_k(r)$ maps $V_i(r)$ into $V_k(r)$.
The mapping $P^i_k(r)\!:\,V_i(r)\to V_k(r)$ is bijective since it is 
invertible. According to \mythetag{5.20}, the mapping 
$P^k_i(r)\!:\,V_k(r)\to V_i(r)$ is inverse to it.\par
     The equality $\dim V_i(r)=\dim V_k(r)$ is an immediate corollary
of the theorem~\mythetheorem{5.3}, which is proved just above. Using
this equality, from the expansions \mythetag{5.16} and \mythetag{5.18}
we derive
$$
\xalignat 2
&\dim V(r)=k_r\,n_r,
&&\dim V(r)=n_r\,\dim V_i(r).
\endxalignat
$$
Comparing these two formulas, we find that
$$
\dim V_i(r)=k_r,\quad i=1,\,\ldots,\,n_r.
$$
The subspaces $V_1,\,\ldots,\,V_{n_r}$ are connected by virtue of the
mappings $P^i_k(r)$. It is easy to verify that the diagram
$$
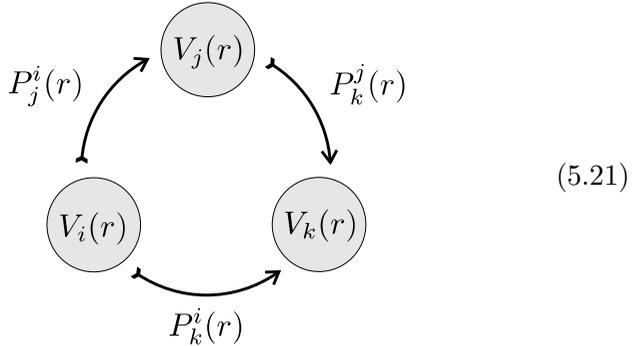
\vcenter to 120pt{\hsize=1cm}
\vadjust{\vskip 5pt\hbox to 0pt{\kern 60pt
\includegraphics{diagram.eps}\hss}\vskip -5pt}
\mytag{5.21}
$$
is commutative. Indeed, $P^j_k(r)\compos P^i_j(r)=P^i_k(r)$. This 
equality is derived from \mythetag{5.11}. Let's choose a basis 
$\bold e^1_1(r),\,\ldots,\,\bold e^1_{k_r}(r)$ in the subspace
$V_1(r)$ and then replicate it to the other subspaces $V_(r)$ by
means of the mappings $P^1_i(r)$: 
$$
\bold e^i_1(r)=P^1_i(r)\bold e^1_1(r),\ .\ .\ .\ ,\ \bold e^i_{k_r}(r)=
P^1_i(r)\bold e^1_{k_r}(r).
$$
Due to the commutativity of the diagram \mythetag{5.21} we get
$$
\hskip -2em
P^i_k\bold e^i_s(r)=\bold e^k_s(r),\quad s=1,\,\ldots,\,k_r.
\mytag{5.22}
$$
The whole set of vectors $\bold e^i_s(r)$, $i=1,\,\ldots,\,n_r$,
$s=1,\,\ldots,\,k_r$ is a linear independent set since the sum of
subspaces in \mythetag{5.18} is a direct sum. Let's define new subspaces
as the spans of the following sets of vectors:
$$
U_s(r)=\langle\bold e^1_s(r),\,\ldots,\,\bold e^{n_r}_s(r)\rangle,
\quad s=1,\,\ldots,\,k_r.\qquad
\mytag{5.23}
$$
Due to the formula \mythetag{5.22} the subspaces \mythetag{5.23} are
invariant with respect to the operators $P^k_i(r)$. However, a stronger
proposition is also valid.
\mytheorem{5.4} The subspace $U_s(r)$ in \mythetag{5.23} is an invariant
subspace of the representation $(\varphi,G,V)$. It is irreducible and the
restriction of $\varphi$ to $U_s(r)$ is isomorphic to $f_r$.
\endproclaim
    In order to prove the invariance of $U_s(r)$ with respect to $\varphi$
we use the relationship \mythetag{5.3}. We have proved it in the very
beginning of this section. Applying it, we get
$$
\hskip -2em
\gathered
\varphi(g)\bold e^i_s(r)=\varphi(g)P^1_i(r)\bold e^1_s(r)=\\
=\sum^{n_r}_{q=1}F^q_i(g,r)\,P^1_q(r)\bold e^1_s(r)=\sum^{n_r}_{q=1}
F^q_i(g,r)\bold e^q_s(r).
\endgathered
\mytag{5.24}
$$
Not only does the relationship \mythetag{5.24} prove the invariance of 
$U_s(r)$ with respect to the operator $\varphi(g)$, but it shows that
the matrix of the operator $\varphi(g)$  in the basis 
$\bold e^1_s(r),\,\ldots,\,\bold e^{n_r}_s(r)$ coincides with the 
matrix of the operator $f_r(g)$. Thus, the second proposition of 
the theorem~\mythetheorem{5.4} is also proved.\par
     As a result we have found the constructive way for expanding 
the space $V$ into a direct sum of subspaces
$$
V=\bigoplus^{m}_{r=1}\bigoplus^{k_r}_{s=1}U_s(r)
$$
that corresponds to the expansion \mythetag{4.14} of the representation
$\varphi$ into its irreducible components.\par
\newpage
\setfirstpage
\topmatter
\title
Contacts
\endtitle
\endtopmatter
\document
\line{\vtop{\hsize 5cm
{\bf Address: }
\medskip\noindent
Ruslan A. Sharipov,\newline
Math. Department,\newline
Bashkir State University,\newline
32 Frunze street,\newline
Ufa 450074, Russia
\medskip
{\bf Phone:}\medskip
\noindent
+7-(347)-273-67-18 (Office)\newline
+7-(917)-476-93-48 (Cell)
}\hss
\vtop{\hsize 4.3cm
{\bf Home address:}\medskip\noindent
Ruslan A. Sharipov,\newline
5 Rabochaya street,\newline
Ufa 450003, Russia
\vskip 1cm
{\bf E-mails:}\medskip
\noindent
r-sharipov\@mail.ru\newline
R\hskip 0.5pt\_\hskip 1.5pt Sharipov\@ic.bashedu.ru
}
}
\bigskip
{\bf URL's:}\medskip
\noindent
\myhref{http://www.geocities.com/r-sharipov/}
{http:/\negskp/www.geocities.com/r-sharipov}\newline
\myhref{http://www.freetextbooks.boom.ru/}
{http:/\negskp/www.freetextbooks.boom.ru}\newline
\myhref{http://sovlit2.narod.ru/}
{http:/\negskp/sovlit2.narod.ru}\newline
\par
\newpage
\setfirstpage
\topmatter
\title
Appendix
\endtitle
\endtopmatter
\document
\rightheadtext{List of publications.}
\leftheadtext{List of publications.}

\enddocument
\end